%% file: EG_of_Curves_upto_two_singular_points_chern_classes.tex
\documentclass[11pt,lettersize,reqno,fullpage]{article}

\usepackage{amsmath,amssymb,amsfonts,mathrsfs,verbatim,enumitem,pstricks,amsthm}
\usepackage[all]{xy}
\usepackage{centernot, xr, accents}
\usepackage{hyperref}
 
\title{{\LARGE Counting curves on a general linear system with up to two singular points}}
\author{Somnath Basu and Ritwik Mukherjee }
\date{}
\usepackage{hyperref}
\hypersetup{
	colorlinks,
	citecolor=black,
	filecolor=black,
	linkcolor=blue,
	urlcolor=black
}

\theoremstyle{plain}
\newtheorem{thm}{Theorem}[section]

\newtheorem{prp}[thm]{Proposition}
\newtheorem{que}[thm]{Question}

\newtheorem{defn}[thm]{Definition}

\newtheorem{rem}[thm]{Remark} 

\newtheorem{mresult}[thm]{Main Theorem}

\newtheorem*{ack}{Acknowledgements}

\include{New_Macros_by_Ritwik_and_Somnath}

\begin{document}

\maketitle

\begin{abstract}
In this paper we obtain an explicit formula for the number of curves
in a compact complex surface $\X$ (passing through the right number of generic points),
that has up to one node and one singularity of codimension $k$, provided 
the total codimension is at most $7$. 
We use a classical fact from 
differential topology: 
the number of zeros of a generic smooth section of a vector 
bundle $V$ over $M$, counted with signs, is 
the Euler class of $V$ evaluated on the fundamental class of $M$.
\end{abstract}

\tableofcontents

\section{Introduction}
\hf\hf Enumeration of singular curves in $\P^2$ (complex projective space) is a classical 
problem in algebraic geometry.
A more general class of enumerative problem is as follows: 
\begin{que}
\label{main_question}
Let $L \lra \X$ be a holomorphic 
line bundle over a compact complex surface and 
$\mathcal{D} :=\P H^0(\X,L) \approx \P^{\delta_L}$ be the space of 
non zero holomorphic sections up to scaling. What is $\N(\A_1^{\delta} \mathfrak{X})$,
the number of curves in $\X$, that belong to the linear system 
$H^0(\X, L)$, passing through $\delta_L-(k+\delta)$ generic points 
and having $\delta$ distinct nodes and one singularity of type $\mathfrak{X}$, 
where $k$ is the codimension of the singularity?
\footnote{By codimension we mean the number of 
conditions having that particular singularity imposes on the space of curves. 
For example, a node is a codimension one 
singularity, a cusp is a codimesnion two singularity, a taconde is a codimension three 
singularity and so on.} 
\end{que}

The main result of this paper (cf. Theorem \ref{one_pt_chern_class_codim_7} and Theorem \ref{two_pt_chern_class_codim_7}) is as follows: 

\begin{mresult}
If $\delta \leq 1$ and  
$\delta+k \leq 7$, then 
we obtain an explicit formula for 
$\N(\A_1^{\delta} \mathfrak{X})$, 
provided $L \lra \X$ is sufficiently ample.
\end{mresult}

Before giving the formulas for $\N(\A_1^{\delta} \mathfrak{X})$, let us 
make a few definitions. 

\begin{defn}
\label{singularity_defn}
Let $L \lra \X$ be a holomorphic 
line bundle over a complex surface and 
and $f \in H^0(\X, L)$ a holomorphic 
section. 
A point $q \in f^{-1}(0)$ \textsf{is of singularity type} $\A_k$,
$\D_k$, $\E_6$, $\E_7$, $\E_8$ or $\XC_8$ if there exists a coordinate system
$(x,y) :(\U,q) \lra (\C^2,0)$ such that $f^{-1}(0) \cap \U$ is given by 
\begin{align*}
\A_k: y^2 + x^{k+1}   &=0  \qquad k \geq 0, \qquad \D_k: y^2 x + x^{k-1} =0  \qquad k \geq 4, \\
\E_6: y^3+x^4 &=0,  \qquad \E_7: y^3+ y x^3=0, 
\qquad \E_8: y^3 + x^5=0, \\
\XC_8: x^4 + y^4 &=0.   
\end{align*}
\end{defn}

\begin{defn}
A holomorphic line bundle $L \lra \X$ over a compact complex manifold $\X$ 
is \textsf{sufficiently $k$-ample} if $L \approx L_1^{\otimes n} \lra \X$, 
where $L_1 \lra \X$ is a very ample line bundle and $n \geq k$. 
\end{defn}

\ni The following two theorems are the main results of this paper:  

\begin{thm}
\label{one_pt_chern_class_codim_7}
Let $\X$ be a two dimensional 
compact complex manifold and $L \lra \X$ a holomorphic line bundle. 
Let 
\[ \mathcal{D} :=\P H^0(\X,L) \approx \P^{\delta_L}, 
\qquad c_1:= c_1(L) \qquad \textnormal{and} \qquad  x_i:= c_i(T^*\X) \] 
where $c_i$ denotes the $i^{th}$ Chern class. 
Denote $\N(\mathfrak{X})$  to be the 
number of curves in $\X$, that belong to the linear system 
$H^0(\X, L)$, passing through $\delta_L-k$ generic 
points and having a singularity of type $\mathfrak{X}$, where $k$ is the codimension of the 
singularity $\mathfrak{X}$. Then 
\begin{align*}
\N(\A_1) &= 3 c_1^2 + 2 c_1 x_1 + x_2,\\ 
\N(\A_2) & = 12 c_1^2 + 12 c_1 x_1 + 2 x_1^2 + 2 x_2,  
\quad \N(\A_3) = 50 c_1^2 + 64 c_1 x_1 + 17 x_1^2 + 5 x_2, \\ 
\N(\A_4) &= 180 c_1^2 + 280 c_1 x_1 + 100 x_1^2, \quad
\N(\A_5) = 630 c_1^2 + 1140 c_1 x_1 + 498 x_1^2 -60 x_2, \\
\N(\A_6) &= 2212 c_1^2 + 4515 c_1 x_1 + 2289 x_1^2 -406 x_2, 
\quad \N(\A_7) = 7812 c_1^2 + 17600 c_1 x_1 + 10022 x_1^2-2058 x_2, \\
\N(\D_4) &= 15 c_1^2 + 20 c_1 x_1 + 5 x_1^2 + 5 x_2, \quad \N(\D_5) = 84 c_1^2 + 132 c_1 x_1 + 44 x_1^2 + 20 x_2, \\
\N(\D_6) &= 224 c_1^2 + 406 c_1 x_1 + 168 x_1^2 + 28 x_2, 
\quad \N(\D_7) = 720 c_1^2 + 1472 c_1 x_1 + 720 x_1^2, \\
\N(\E_6) &= 84 c_1^2 + 147 c_1 x_1 + 57 x_1^2 + 18 x_2, 
\quad \N(\E_7) = 252 c_1^2 + 488 c_1 x_1 + 217 x_1^2 + 42 x_2,
\end{align*}
provided $L$ is sufficiently $\ct_{\mathfrak{X}}$-ample, where 
\bgd
\ct_{\A_k} := k+1, ~~\ct_{\D_k} := k-1, ~~\ct_{\E_6} = 4, ~~\ct_{\E_7} = 4.
\edd
\end{thm}
\begin{thm}
\label{two_pt_chern_class_codim_7}
Let $\X$ be a two dimensional 
compact complex manifold and $L \lra \X$ a holomorphic line bundle. 
Let 
\[ \mathcal{D} :=\P H^0(\X,L) \approx \P^{\delta_L}, 
\qquad c_1:= c_1(L) \qquad \textnormal{and} \qquad  x_i:= c_i(T^*\X) \] 
where $c_i$ denotes the $i^{th}$ Chern class. 
Denote $\N(\A_1\mathfrak{X})$ to be the 
number of curves in $\X$, that belong to the linear system 
$H^0(\X, L)$, passing through $\delta_L-(k+1)$ generic 
points and having one simple node and one singularity of type $\mathfrak{X}$, where $k$ is the codimension of the 
singularity $\mathfrak{X}$. 
Then 
\begin{align*}
\N(\A_1\A_1) &= -42 c_1^2 + 9 c_1^4 -39 c_1 x_1 + 12 c_1^3 x_1 -6 x_1^2 + 4 c_1^2 x_1^2 - 7 x_2 + 6 c_1^2 x_2 + 
4 c_1 x_1 x_2 + x_2^2, \\
\N(\A_1 \A_2) &= -240 c_1^2 + 36 c_1^4 -288 c_1 x_1 + 60 c_1^3 x_1 -72 x_1^2 + 30 c_1^2 x_1^2+ 4 c_1 x_1^3 -24 x_2
                 + 18 c_1^2 x_2\\ 
              &~~ + 16 c_1 x_1 x_2 + 2 x_1^2 x_2 + 2 x_2^2,\\
\N(\A_1\A_3) &= -1260 c_1^2 + 150 c_1^4 -1820 c_1 x_1 + 292 c_1^3 x_1 -596 x_1^2 + 179 c_1^2 x_1^2 + 34 c_1 x_1^3
-60 x_2\\ 
             &  ~~ + 65 c_1^2 x_2 + 74 c_1 x_1 x_2 + 17 x_1^2 x_2 + 5 x_2^2, \\
\N(\A_1\A_4) &= -5460 c_1^2 + 540 c_1^4 -9240 c_1 x_1 + 1200 c_1^3 x_1 -  3740 x_1^2 + 860 c_1^2 x_1^2 + 200 c_1 x_1^3 
+200 x_2\\ 
             & ~~ + 180 c_1^2 x_2 + 280 c_1 x_1 x_2 + 100 x_1^2 x_2, \\
\N(\A_1\A_5) &= -22428 c_1^2 + 1890 c_1^4 -43197 c_1 x_1 + 4680 c_1^3 x_1 -20535 x_1^2 + 3774 c_1^2 x_1^2+ 996 c_1 x_1^3\\ 
             & ~~ +2754 x_2 + 450 c_1^2 x_2 + 1020 c_1 x_1 x_2 + 498 x_1^2 x_2 -60 x_2^2, \\
\N(\A_1\A_6) &= -90468 c_1^2 + 6636 c_1^4 -193816 c_1 x_1 + 17969 c_1^3 x_1 -104503 x_1^2 + 15897 c_1^2 x_1^2 \\ 
             &~~ +4578c_1 x_1^3 + 18522 x_2 + 994 c_1^2 x_2 + 3703 c_1 x_1 x_2 + 2289 x_1^2 x_2 -406 x_2^2, \\
\N(\A_1\D_4) &= -420 c_1^2 + 45 c_1^4 -624 c_1 x_1 + 90 c_1^3 x_1 -196 x_1^2 + 55 c_1^2 x_1^2 + 10 c_1 x_1^3\\ 
             &~~ -100x_2 + 30 c_1^2 x_2 + 30 c_1 x_1 x_2 + 5 x_1^2 x_2 + 5 x_2^2, \\
\N(\A_1\D_5) &= -2688 c_1^2 + 252 c_1^4 -4564 c_1 x_1 + 564 c_1^3 x_1 -1744 x_1^2 + 396 c_1^2 x_1^2 + 88 c_1 x_1^3 \\ 
             & ~~ -456x_2 + 144 c_1^2 x_2 + 172 c_1 x_1 x_2 + 44 x_1^2 x_2 + 20 x_2^2,\\ 
\N(\A_1 \D_6) &= -8316 c_1^2 + 672 c_1^4 -16008 c_1 x_1 + 1666 c_1^3 x_1 - 7281 x_1^2 + 1316 c_1^2 x_1^2 + 336 c_1 x_1^3 \\ 
              & ~~ -546 x_2 + 308 c_1^2 x_2 + 462 c_1 x_1 x_2 + 168 x_1^2 x_2 + 28 x_2^2, \\
\N(\A_1 \E_6) &= -2916 c_1^2 + 252 c_1^4 -5400 c_1 x_1 + 609 c_1^3 x_1 -2295 x_1^2 + 465 c_1^2 x_1^2 + 114 c_1 x_1^3\\ 
              &~~ -486 x_2 + 138 c_1^2 x_2 + 183 c_1 x_1 x_2 + 57 x_1^2 x_2 + 18 x_2^2,
\end{align*}
provided $L$ is sufficiently $\ct_{\A_1\mathfrak{X}}$-ample, where 
$\ct_{\A_1 \mathfrak{X}} : = 2 +  \ct_{\mathfrak{X}}$.
\end{thm}


\begin{rem}
In the formulas presented in Theorems \ref{one_pt_chern_class_codim_7} and  
\ref{two_pt_chern_class_codim_7}, 
if there is a degree four cohomology class, 
we mean the cohomology class evaluated on the fundamental class $[\X]$. When 
there is a degree eight cohomology class, we mean the cohomology class 
evaluated on $[\X \times \X]$. 
\end{rem}

\section{A description of the method used in this paper} 
\label{method_description}
\hf\hf Let us first recall the results of our earlier papers. In \cite{BM13} and 
\cite{BM13_2pt_published}, we 
prove the special cases of Theorems \ref{one_pt_chern_class_codim_7} and 
\ref{two_pt_chern_class_codim_7} respectively when $\X:= \P^2$ and 
$L := \gamma_{\P^2}^{* d}$. In \cite{RM_Hypersurfaces},  
we obtain the first three formulas presented in Theorem \ref{one_pt_chern_class_codim_7}.
\footnote{Actually the results in \cite{RM_Hypersurfaces} are slightly more general; we obtain an 
enumerative formula for \textit{hypersurfaces} with a node, a cusp or a tacnode.} \\ 
\hf \hf We will now give a description of the method we use to obtain the enumerative formulas 
presented in this paper. Our goal is to enumerate curves in a linear system, passing through the right 
number of generic points, with  certain singularities. 
A curve having a singularity implies that a certain derivative is zero. We interpret this 
derivative as the section of some vector bundle. If  
this section is transverse to the zero set, our 
desired number is the Euler class of the vector bundle. Computing the Euler class is 
completely elementary via the splitting principle.\\ 
\hf \hf However, it turns out that there is a subtlety involved very soon in this process. 
One may look at a simple example to understand what is going on. Consider the complex  
polynomial 
\[ f(z):= z^5(z-1)(z-2)(z-3)\] 
and let us ask what is $n$, i.e., 
\[ n:= \big|\{ z \in \mathbb{C}: f(z) =0, ~~z \neq 0\} \big |=\,? \]
A first guess would be $8$ by looking at the degree of $f$. However, this answer is incorrect 
because $f$ vanishes when $z=0$. We need to compute the contribution of $f$ to its degree 
from the point $z=0$. To compute this contribution, we smoothly perturb the function $f$ and 
count (with a sign) how many zeros are there in a small neighborhood of $z=0$, i.e., 
we count the signed number of solutions of 
\[ f(z) + \nu(z) = 0, \qquad |z-0| < \epsilon \]
where $\nu$ is a small generic perturbation (i.e., the $C^0$-norm of $\nu$ is small) 
and $\epsilon$ is sufficiently small. This number is $5$. Hence 
\[ \textnormal{deg}(f) = n + 5 \qquad \implies \qquad n = 8-5. \]
\hf \hf To return to our main discussion of enumerating curves with a singularity, 
it turns out that once the singularity 
becomes too degenerate, or we allow the curve to have more than one singular point, the Euler 
class counts too much. 
There is a boundary contribution which we have to subtract off from 
the Euler class. This part is highly non-trivial and challenging.  \\ 
\hf \hf In  \cite{Z1}, \cite{BM13} and \cite{BM13_2pt_published} the authors 
carry out this topological method 
and 
obtain the special cases of Theorems \ref{one_pt_chern_class_codim_7} and 
\ref{two_pt_chern_class_codim_7} respectively when $\X:= \P^2$ and 
$L := \gamma_{\P^2}^{* d}$ (when $d$ is sufficiently high). The bound on 
$d$ is required to ensure that the relevant sections are transverse to the 
zero set. 
It may seem that one of the potential 
difficulties of generalizing the formulas to any arbitrary line bundle $L \lra \X$ 
is that we do not know if these sections are going to be transverse to the 
zero set. As it turns out there is a very simple criteria to guarantee transversality; 
if the line bundle is sufficiently ample, then the relevant sections are going to be transverse. 
Once transversality is achieved, the arguments given 
in \cite{BM13} and \cite{BM13_2pt_published} apply immediately to the more general setup 
of considering curves in a linear system. As a result, we obtain these formulas in terms of 
Chern classes.

\section{Organization of this paper} 

\hf\hf We now describe the basic organization of this paper. As explained in section \ref{method_description}, 
the two main aspects of applying our method is: proving transversality 
and computing the boundary contribution to the Euler class. Towards showing transversality, in \cite{BM13}, \cite{BM13_2pt_published} and \cite{BM_Detail}, we give 
 a rigorous proof of why the relevant sections are transverse in the special case when 
 $\X := \P^2 $ and $L := \gamma_{\P^2}^{* d}$ (provided $d$ is sufficiently large). In \cite{RM_Hypersurfaces}, 
 the author shows that in the more general setup, the relevant sections that arise during the 
 computations of $\N(\A_1)$, $\N(\A_2)$ and $\N(\A_3)$ are also going to be transverse provided 
 the line bundle is sufficiently ample. Moreover, the argument that was employed was essentially 
 mimicking the arguments employed in the case of projective space.  In particular, 
 it is easy to see that once the line bundle is sufficiently ample, all the arguments presented in  
\cite{BM_Detail} to prove transversality can be mimicked for the more general setup of $L \lra \X$. 
Hence, we have decided to omit the proof of why the relevant sections are transverse from this paper. 
The reader can can refer to  \cite{BM_Detail} and \cite{RM_Hypersurfaces} to see why transversality 
holds. \\ 
\hf \hf We next discuss the more crucial aspect of computing the boundary contribution to the 
Euler class. As we explained in section \ref{method_description}, the boundary contribution 
was computed rigorously in \cite{BM13} and \cite{BM13_2pt_published} for the special 
case of the projective space. A little bit of thought 
shows that as long as the relevant sections are transverse to the zero set, 
exactly the same arguments apply to the more general setup of $L \lra \X$ to compute the 
boundary contribution. Hence, we omit the proof of how we obtain the multiplicities from each 
boundary component; we explicitly state the multiplicities and show how to obtain the final 
formula.  The reader can can refer to  \cite{BM13} and \cite{BM13_2pt_published} to see 
how those multiplicities were actually obtained. 

\begin{ack}
The second author is grateful to Aleksey Zinger for pointing out \cite{Z1} and 
explaining the topological method employed in that paper.
\end{ack}

\section{A survey of related results in Enumerative Geometry} 

\hf\hf We now give a  brief survey of related results in this area of mathematics, 
namely Enumerative Geometry 
of Singular Curves and Hypersurfaces. 
We start by looking at the results of M.Kazarian. We should mention at the outset that 
although we are not completely certain, 
we believe it is very likely that by applying the methods described in \cite{Kaz}, 
the results of
Theorem \ref{one_pt_chern_class_codim_7} 
and Theorem \ref{two_pt_chern_class_codim_7} can be recovered. 
However, Kazarian's methods 
\textit{are completely different from ours}. 
Furthermore, we believe that our method complements his 
method  very well, since each method has its own  
advantages and disadvantages. \\ 
\hf \hf Let us now explain the method of M.Kazarian. 
His method works on the principle that there exists a universal formula for these 
enumerative numbers in terms of the Chern classes. He then goes on to consider 
enough special cases to find out what that exact combination is.\footnote{To take a simple 
example; suppose there is a polynomial of degree $m$. To find out what the polynomial is, 
we simply have to find the value of the polynomial at enough points.}  
One of the difficulties of this method is to prove the existence of such a universal 
formula. We also believe it is usually difficult  
to think of enough special cases in a 
given situation. However, this method has been successfully applied  
in many situations and in particular we believe it recovers our results. \\
\hf\hf The reader who has read section \ref{method_description} will see immediately 
that this method is completely different from ours; \textit{we do not make} 
\textit{any assumption} that there is a universal formula in terms of Chern classes.
Aside from the results of Kazarian, the rest of the results in this field are either 
\textit{special cases} of Theorems
\ref{one_pt_chern_class_codim_7} and \ref{two_pt_chern_class_codim_7}  
or a \textit{completely different class} of enumerative problems. \\
\hf \hf Let us look at the results of Dmitry Kerner. In his paper \cite{Ker1}, Kerner 
considers the special case of 
Theorem \ref{one_pt_chern_class_codim_7}, when 
$\X:= \P^2$ and $L := \gamma_{\P^2}^{* d}$. 
Our results are consistent with his    
in this special case.
Kerner also considers 
in his paper \cite{Ker2} 
the special cases of Theorem \ref{two_pt_chern_class_codim_7}, when 
$\X:= \P^2$ and $L := \gamma_{\P^m}^{* d}$ and obtains three of the formulas we have stated 
in Theorem \ref{two_pt_chern_class_codim_7}; a formula for $\N(\A_1 \A_1)$, 
$\N(\A_1\A_2)$ and $\N(\A_1\D_4)$. \\
\hf \hf The crucial difference between the results of Kazarian and 
our results and those  
obtained in \cite{Ker1}, \cite{Ker_Hypersurface} and \cite{Ker2} is that the author
there obtains results \textit{only} for the special case of the linear system 
$\gamma_{\P^2}^{*d} \lra \P^2$, while Kazarian and our results are for \textit{any} linear 
system $L \lra \X$ that is sufficiently ample. \\
\hf \hf Next, let us look at the results of I.Vainsencher. 
He considers a different class of enumerative problems 
(with the exception of $\N(\A_1)$ and $\N(\A_2)$;
our results are consistent with his).
In his paper \cite{Van}, Vainsencher 
considers a general linear system $L \lra \X$ and enumerates curves that have up to six nodes. 
He also obtains a formula for $\N(\A_2)$ in his paper \cite{Van_Cusp}.\\ 
\hf \hf Next, let us look at the results  
S.Klienman and R.Piene. They do obtain few of the 
formulas we have obtained in Theorems \ref{one_pt_chern_class_codim_7} 
and \ref{two_pt_chern_class_codim_7}; namely $\N(\A_1), \N(\D_4), \N(\D_6), \N(E_7), \N(\A_1\A_1)$, 
$\N(\A_1\D_4)$ and $\N(\A_1 \D_6)$. Our 
results are consistent with theirs. 
They also study a different class of enumerative questions, namely 
enumerating curves that have up to eight simple nodes, or one triple point 
and up to four simple nodes, or one $\D_6$ node and up to two simple nodes. 
Our results are of a different nature from theirs; in Theorems 
\ref{one_pt_chern_class_codim_7}, \ref{two_pt_chern_class_codim_7} we enumerate curves 
with up to one node and one arbitrarily degenerate singularity (till a total codimension of seven). \\
\hf \hf Finally, we note that in their papers 
\cite{Ran1}, \cite{Ran2} and \cite{CH}, Z. Ran, L.Caporasso and J.Harris have 
obtained a formula for the number of degree $d$-curves in $\P^2$ (through the right 
number of generic points) having $\delta$-nodes, for any $\delta$. However, their 
results are only for $\P^2$. Moreover, the allowed singularities 
in their cases are simple nodes and not anything more degenerate. \\
\hf \hf To summarize, aside from the results of Kazarian, all the other results are either 
\textit{special 
cases} of Theorems 
\ref{one_pt_chern_class_codim_7} and \ref{two_pt_chern_class_codim_7} or are 
enumerative results of a different nature. 
Our method is completely different from that of Kazarian 
and complements it very well. Moreover, 
our method has the potential to go much further beyond codimension seven 
(just like the method of Kazarian has the potential to go a lot further).  

\section{Topological computations: one singular point} 
\label{top_comp_one_sing_pt}
\hf\hf In this section we will give a proof of Theorem \ref{one_pt_chern_class_codim_7}. 
Let us first set up some  notation. Given a singularity $\mf{X}$ let us define the 
following spaces 
\begin{align*}
\mathfrak{X} &:= \{ ([f], q) \in \D \times \X: f ~~\textnormal{has a signularity of type $\mathfrak{X}$ at $q$} \}, 
\\
\hat{\mathfrak{X}} &:= \{ ([f], l_q) \in \D \times \mathbb{P}T\X: ([f], q) \in \mathfrak{X} \}, 
\\
\mathcal{P}\A_k &:= \{ ([f], l_q) \in \D \times \mathbb{P} T\X: ([f], q) \in \A_k, ~~\nabla^2 f(v, \cdot) =0 
\qquad \forall v \in l_q\} \qquad \textnormal{if ~~$k \geq 2$}, \\
\mathcal{P}\D_4 &:= \{ ([f], l_q) \in \D \times \mathbb{P} T\X: ([f], q) \in \D_4, ~~\nabla^3 f(v,v, v) =0 
\quad \forall v \in l_q\}, \\
\mathcal{P}\D_k &:= \{ ([f], l_q) \in \D \times \mathbb{P} T\X: ([f], q) \in \D_k, ~~\nabla^3 f(v,v, \cdot) =0 
\quad \forall v \in l_q\} \qquad \textnormal{if ~~$k \geq 5$}, \\
\mathcal{P}\E_k &:= \{ ([f], l_q) \in \D \times \mathbb{P} T\X: ([f], q) \in \E_k, ~~\nabla^3 f(v,v, \cdot) =0 
\quad \forall v \in l_q\} \qquad \textnormal{if ~~$k=6,7,8$}.
\end{align*}

Next, if $\mf{X}$ is a codimension $k$ singularity, 
we define the following two numbers 
\begin{align*}
\N(\mf{X}, n_1, m_1, m_2) &:= 
\lan c_1^{n_1} x_1^{m_1} x_2^{m_2} y^{\delta_L-(n_1+m_1+2m_2 + k)},  
~[\overline{\mf{X}}]\ran,  \\
\N(\mathcal{P} \mf{X}, n_1, m_1, m_2, \theta) &:= 
\lan c_1^{n_1} x_1^{m_1} x_2^{m_2} \lambda^{\theta}y
^{\delta_L-(n_1+m_1+2m_2 +\theta+ k)},  
~[\overline{\mc{P} \mf{X}}]\ran, 
\end{align*}
where 
\begin{align*}
c_1 &: = c_1(L), ~~x_i := c_i(T^*\X), ~~\lambda := c_1(\hat{\gamma}^*), ~~ y:= c_1(\gamma_{\D}^*) 
\end{align*}
and $\gamma_{\D} \lra \D$ and $\hat{\gamma} \lra \mathbb{P} T\X$ are 
the tautological line bundles. \\
\hf \hf We will now give a series of formulas to compute $\N(\A_1, n_1, m_1, m_2)$ and 
$\N(\mc{P}\mf{X}, n_1, m_1, m_2, \theta)$. 
Note that $\N(\A_1) = \N(\A_1, 0,0,0)$.
In order to compute the remaining $\N(\mf{X})$ we 
do the following: if $\mf{X}$ is anything other than $\D_4$, then we observe that 
$\N(\mf{X}) = \N(\mc{P}\mf{X}, 0, 0, 0, 0)$. If $\mf{X} = \D_4$, then we observe that 
$\N(\D_4) = \frac{\N(\mc{P}\D_4,0,0,0)}{3}$.\\  
\hf \hf Note that, Propositions  \ref{a1_chern} to \ref{pa7_chern} 
and Propositions \ref{a1a1_chern} to \ref{a1pd6_chern}
give \textsf{recursive} formulas 
to compute $\N(\mf{X})$ and $\N(\A_1\mf{X})$. 
We wrote a Mathematica program 
to obtain the final formulas given in Theorem 
\ref{one_pt_chern_class_codim_7} and \ref{two_pt_chern_class_codim_7}. 
The program can be obtained from our homepage 
\[ \textnormal{\url{https://www.sites.google.com/site/ritwik371/home}}. \]
\hf \hf Before we prove Theorem \ref{one_pt_chern_class_codim_7},  
we also need to define  the following spaces; they will come up during the course of 
our computations: 
\begin{align*}
\hat{\A}_1^{\#}  := \{ ([f], l_q) \in \D \times \P T\X &: 
 f(q) =0, \nabla f|_q =0,  \nabla^2 f|_q(v, \cdot) \neq 0, \forall ~v \in l_{q}-0  \}  \\
\hat{\D}_4^{\#}  := \{ ([f], q) \in \D \times \P T\X &: 
f(q) =0, \nabla f|_q =0, \nabla^2 f|_q \equiv 0, \nabla^3 f|_q (v,v,v) \neq 0, \forall ~v \in l_{q}-0 \} \\
\hat{\D}_k^{\#\flat}  := \{ ([f], l_q) \in \D \times \P T\X &: \textnormal{$f$ has a $\D_k$ singularity at $q$},  
~~\nabla^3 f|_q (v,v,v) \neq 0, \forall ~v \in l_{q}-0,\,\,k\geq 4\}. \\
\hat{\XC}_8^{\#} := \{ ([f], l_q) \in \D \times \P T\X &: 
f(q) =0, \nabla f|_q =0, \nabla^2 f|_q \equiv 0, \nabla^3 f|_q =0,  \\ 
                & \qquad \nabla^4 f|_q (v,v,v, v) \neq 0 ~\forall ~v \in l_{q}-0 \}, \\   
\PP\D_k^{\vee} := \{ ([f], q) \in \D \times \P T\X &: ([f], q) \in \D_k, 
~~\nabla^3 f|_q(v,v,v) =0, ~~\nabla^3 f|_q(v,v,w) \neq 0,  \\
&  \forall ~~v  \in  l_{q} -0 ~~\textup{and} ~~w \in (T_{q}\X)/l_{q} -0 \}, \qquad \textnormal{if $~k>4$}. 
\end{align*}
\hf \hf We will use the following 
fact from differential topology 
(cf. \cite{BoTu}, Proposition 12.8):
\begin{thm} 
\label{Main_Theorem} 
Let $V\lra M$ be an oriented  vector bundle over a compact oriented manifold $M$ 
and $s:M \lra V$ a smooth section that is transverse to the zero set. 
Then the 
Poincar\'{e} dual of $[s^{-1}(0)]$ in $M$ is the Euler class of $V$. 
In particular, if the rank of $V$ is same as the dimension of $M$, 
then the signed cardinality of $s^{-1}(0)$ is the Euler class of $V$, evaluated on the fundamental class of 
$M$.
\end{thm}

\ni We are now ready to prove Theorem \ref{one_pt_chern_class_codim_7}. 
It is to be understood that  Propositions  \ref{a1_chern} to \ref{pa7_chern} 
and Propositions \ref{a1a1_chern} to \ref{a1pd6_chern} are true provided 
$L$ is appropriately ample (as stated in Theorems \ref{one_pt_chern_class_codim_7} 
and \ref{two_pt_chern_class_codim_7}).


\begin{prp}
\label{a1_chern}
The number $\N(\A_1, n_1, m_1, m_2)$ is given by 
\begin{align}
\label{a1_chern_eqn}
\N(\A_1, n_1, m_1, m_2) & = \begin{cases} 3 c_1^2 + 2 c_1 x_1 + x_2 & \mbox{if} ~~(n_1, m_1, m_2) =(0,0,0), \\ 
   3c_1^2 + c_1 x_1  & \mbox{if }  ~~(n_1, m_1, m_2) = (1,0,0),\\ 
   c_1^2 & \mbox{if } ~~(n_1, m_1, m_2) = (2,0,0), \\ 
   3 c_1 x_1 + x_1^2 &  \mbox{if} ~~~(n_1, m_1, m_2) = (0,1,0), \\ 
   c_1 x_1 & \mbox{if} ~~~(n_1, m_1, m_2) = (1,1,0), \\ 
   x_1^2 & \mbox{if} ~~~(n_1, m_1, m_2) = (0,2,0), \\ 
   x_2 & \mbox{if} ~~~(n_1, m_1, m_2) = (0,0,1), \\ 
   0 & \mbox{otherwise}.\end{cases} 
\end{align}
\end{prp}
\pf Let $\mu$ be a generic pseudocycle representing the homology class 
Poincar\'{e} dual 
to 
\[c_1^{n_1} x_1^{m_1} x_2^{m_2} y^{\delta_L - (n_1 + m_1 + 2 m_2 + 1)}.\] 
We now define sections of the following two bundles 
\begin{align}
\psi_{\A_0}:  (\D \times \X) \cap \mu \lra \mathcal{L}_{\A_0} & := \gamma_{\D}^* \otimes L \lra \D \times \X, 
\qquad \textnormal{given by} \qquad  
\{\psi_{\A_0}([f], q)\}(f) := f(q) \qquad \textnormal{and} \nonumber \\
\psi_{\A_1} : \psi_{\A_0}^{-1}(0) \lra \mathcal{V}_{\A_1} &:= \gamma_{\D}^*\otimes T^*\X \otimes L, 
\qquad \textnormal{given by} \qquad 
\{\psi_{\A_1}([f], q)\}(f) := \nabla f|_q. \label{psi_a0_a1_section_defn}
\end{align}
In \cite{RM_Hypersurfaces} we show that if $L$ is sufficiently $2$-ample, then 
these sections are transverse to the zero set. Hence 
\begin{align}
\N(\A_1, n_1, m_1, m_2) & = \lan e(\mathcal{L}_{\A_0}) e(\mathcal{V}_{\A_1}), ~~[\D \times \X] \cap [\mu] \ran.  
\label{na1_pseudocycle}
\end{align}
Equation \eqref{na1_pseudocycle} and the 
Splitting Principle, imply \eqref{a1_chern_eqn}. \qed

\begin{prp}
The number $\N(\mc{P}\A_2, n_1, m_1, m_2, \theta)$ is given by
\begin{align}
\label{pa2_chern_eqn}
\N(\mc{P}\A_2, n_1, m_1, m_2, \theta) &= 
\begin{cases} 2\N(\A_1, n_1, m_1, m_2) + 2 \N(\A_1, n_1, m_1+1, m_2)\\ 
+ 2\N(\A_1, n_1+1, m_1, m_2) \qquad \mbox{if} ~~\theta =0,\\ \\
                                  \N(\A_1, n_1, m_1, m_2) + 2\N(\A_1, n_1+1, m_1, m_2) + \N(\A_1, n_1+2, m_1, m_2) \\ 
                                  + 3 \N(\A_1, n_1, m_1+1, m_2) + 3 \N(\A_1, n_1+1, m_1+1, m_2) \\ 
                                  + 
                                  2\N(\A_1, n_1, m_1+2, m_2) 
                                  \qquad \mbox{if} ~~\theta =1,\\ \\
                                  \N(\mc{P}\A_2, n_1, m_1+1, m_2, \theta-1) \\
                                 -\N(\mc{P}\A_2, n_1, m_1, m_2+1, \theta-2)
                                 \qquad \mbox{if} ~~\theta >1. \end{cases}
\end{align}
\end{prp}
\pf Let $\mu$ be a generic pseudocycle representing the homology class 
Poincar\'{e} dual 
to 
\[c_1^{n_1} x_1^{m_1} x_2^{m_2} \lambda^{\theta}y^{\delta_L - (n_1 + m_1 + 2 m_2 + 2)}.\] 
We now define a section of the following bundle 
\begin{align*}
\us_{\PP \A_2}: & \ov{\hat{\A}}_1 \cap \mu \lra  \mathbb{V}_{\mc{P}\A_2} := 
\hat{\gamma}^* \otimes \gamma_{\D}^* \otimes T^*\X \otimes L, \qquad \textnormal{given by} \\ 
\{ \us_{\PP \A_2}([f], l_q) \}(v \otimes f) &:= \nabla ^2 f(v, \cdot) \qquad \forall ~~v \in l_q.
\end{align*}
In \cite{RM_Hypersurfaces} we show that if $L$ is sufficiently $3$-ample, then 
this section is transverse to the zero set. Hence 
\begin{align}
\N(\PP\A_2, n_1, m_1, m_2, \theta) & = 
\Big\lan e(\mathbb{V}_{\PP \A_2}), ~~[\ov{\hat{\A}}_1] \cap [\mu] \Big\ran.  
\label{npa2_pseudocycle}
\end{align}
Equation \eqref{npa2_pseudocycle}, the Splitting Principle and the ring structure 
of $H^*(\P T\X)$ imply \eqref{pa2_chern_eqn}. \qed 

\begin{prp}
The number $\N(\mc{P}\A_3, n_1, m_1, m_2, \theta)$ is given by 
\begin{align}
\N(\mc{P}\A_3, n_1, m_1, m_2, \theta) &= 3\N(\mc{P}\A_2, n_1, m_1, m_2, \theta+1) 
+\N(\mc{P}\A_2, n_1, m_1, m_2, \theta) \nonumber \\  
& ~~ + \N(\mc{P}\A_2, n_1+1, m_1, m_2, \theta).  \label{pa3_chern_eqn}
\end{align}
\end{prp}
\pf Let $\mu$ be a generic pseudocycle representing the homology class 
Poincar\'{e} dual 
to 
\[c_1^{n_1} x_1^{m_1} x_2^{m_2} \lambda^{\theta}y^{\delta_L - (n_1 + m_1 + 2 m_2 + 3)}.\] 
We now define a section of the following bundle 
\begin{align*}
\us_{\PP \A_3}: & \ov{\PP \A}_2 \cap \mu \lra  \mathbb{L}_{\mc{P}\A_3} := 
\hat{\gamma}^{*3} \otimes \gamma_{\D}^* \otimes L, \qquad \textnormal{given by} \\ 
\{ \us_{\PP \A_3}([f], l_q) \}(v^{\otimes 3} \otimes f) &:= \nabla ^3 f(v, v, v) \qquad \forall ~~v \in l_q.
\end{align*}
In \cite{RM_Hypersurfaces} we show that if $L$ is sufficiently $4$-ample, then 
this section is transverse to the zero set. Hence 
\begin{align}
\N(\PP\A_3, n_1, m_1, m_2, \theta) & = 
\big\lan e(\mathbb{L}_{\PP \A_3}), ~~[\ov{\PP \A}_2] \cap [\mu] \big\ran,  
\label{npa3_pseudocycle}
\end{align}
which gives us \eqref{pa3_chern_eqn}. \qed

\begin{prp}
The number $\N(\mc{P}\D_4, n_1, m_1, m_2, \theta)$ is given by 
\begin{align}
\N(\mc{P}\D_4, n_1, m_1, m_2, \theta) &= 2 \N(\mc{P}\A_3, n_1, m_1+1, m_2, \theta) 
-2 \N(\mc{P}\A_3, n_1, m_1, m_2, \theta+1) \nonumber \\  
& ~~ + \N(\mc{P}\A_3, n_1, m_1, m_2, \theta) + \N(\mc{P} \A_3, n_1+1, m_1, m_2, \theta).  \label{pd4_chern_eqn}
\end{align}
\end{prp}
\pf Let $\mu$ be a generic pseudocycle representing the homology class Poincar\'{e} dual 
to 
\[c_1^{n_1} x_1^{m_1} x_2^{m_2} \lambda^{\theta} y^{\delta_L - (n_1 + m_1 + 2 m_2 + \theta+4)}.\] 
It is shown in \cite{BM13} that 
\begin{align*}
\ov{\PP \A}_3 &= \PP \A_3 \cup \ov{\PP \A}_4 \cup \ov{\mc{P}\D}_4.
\end{align*}
We now define a section of the following bundle 
\begin{align}
\Psi_{\mathcal{P} \D_4}:  \ov{\mc{P}\A}_3 \cap \mu \lra \mathbb{L}_{\mathcal{P} \D_4} & := 
(T\X/\hat{\gamma})^{*2} \otimes \gamma_{\D}^* \otimes L, \qquad \textnormal{given by} \nonumber \\ 
\{\Psi_{\mathcal{P} \D_4}([f], l_{q})\}(w^{\otimes 2} \otimes f) & := \nabla^2 f|_{q}(w,w) \qquad 
\forall ~w \in T_q\X/l_q.
\label{psi_pd4_section_defn}
\end{align}
If $L$ is sufficiently $3$-ample, then this section  
vanishes on $\PP\D_4 \cap \mu$ transversally. Moreover, it does not vanish on 
$\PP \A_4  \cap \mu$. Since $\mu$ is generic, we conclude that 
$\ov{\PP \A}_4  \cap \mu = \PP \A_4  \cap \mu$; hence this section does not vanish 
on $\ov{\PP \A}_4  \cap \mu$. Hence, 
\begin{align*}
\N(\mc{P}\D_4, n_1, m_1, m_2, \theta)  &= \big\lan e(\mathbb{L}_{\mathcal{P} \D_4}), ~~[\ov{\mc{P}\A}_3] \cap [\mu] \big\ran.  
\end{align*}
This proves \eqref{pd4_chern_eqn}. \qed 

\begin{prp}
The number $\N(\mc{P}\D_5, n_1, m_1, m_2, \theta)$   is given by 
\begin{align}
\N(\mc{P}\D_5, n_1, m_1, m_2, \theta) &=  \N(\mc{P}\D_4, n_1, m_1, m_2, \theta+1) 
+\N(\mc{P}\D_4, n_1, m_1+1, m_2, \theta) \nonumber \\  
& ~~ + \N(\mc{P}\D_4, n_1, m_1, m_2, \theta) + \N(\mc{P} \D_4, n_1+1, m_1, m_2, \theta).  \label{pd5_chern_eqn}
\end{align}
\end{prp}
\pf Let $\mu$ be a generic pseudocycle representing the homology class Poincar\'{e} dual 
to 
\[c_1^{n_1} x_1^{m_1} x_2^{m_2} \lambda^{\theta} y^{\delta_L - (n_1 + m_1 + 2 m_2 + \theta+5)}.\] 
It is shown in \cite{BM13} that 
\begin{align*}
\ov{\mc{P}\D}_4 &= \PP \D_4 \cup \ov{\mc{P}\D}_5 \cup \ov{\mc{P}\D^{\vee}_5}.  
\end{align*}
We now define a section of the following bundle 
\begin{align}
\Psi_{\mathcal{P} \D_5}:  \ov{\mc{P}\D}_4 \cap \mu \lra \mathbb{L}_{\mathcal{P} \D_5} & := 
\hat{\gamma}^{* 2} \otimes (T\X/\hat{\gamma})^{*} \otimes \gamma_{\D}^* \otimes L, \qquad 
\textnormal{given by} \nonumber \\ 
\{\Psi_{\mathcal{P} \D_5}([f], l_{q})\}(v^{\otimes 2} \otimes w \otimes f) & := \nabla^3 f|_{q}(v,v,w) \qquad 
\forall ~~v \in l_q, ~~w \in T_q\X/l_q.
\label{psi_pd4_section_defn}
\end{align}
If $L$ is sufficiently $4$-ample, then the section $\Psi_{\mathcal{P} \D_5}$ 
restricted to $\PP \D_5 \cap \mu$ vanishes transversally. Moreover, it does not vanish 
on $\PP \D_5^{\vee} \cap \mu$. 
Hence 
\begin{align*}
\N(\mc{P}\D_5, n_1, m_1, m_2, \theta)  &= \big\lan e(\mathbb{L}_{\mathcal{P} \D_5}), ~~[\ov{\mc{P}\D}_4]\cap [\mu] \big\ran.  
\end{align*}
This proves \eqref{pd5_chern_eqn}. \qed 

\begin{prp}
The number $\N(\mc{P}\E_6, n_1, m_1, m_2, \theta)$   is given by 
\begin{align}
\N(\mc{P}\E_6, n_1, m_1, m_2, \theta) &=  2\N(\mc{P}\D_5, n_1, m_1+1, m_2, \theta) 
-\N(\mc{P}\D_5, n_1, m_1, m_2, \theta+1) \nonumber \\  
& ~~ + \N(\mc{P}\D_5, n_1+1, m_1, m_2, \theta)+ \N(\mc{P}\D_5, n_1, m_1, m_2, \theta).  \label{pe6_chern_eqn}
\end{align}
\end{prp}
\pf Let $\mu$ be a generic pseudocycle representing the homology class Poincar\'{e} dual 
to 
\[c_1^{n_1} x_1^{m_1} x_2^{m_2} \lambda^{\theta} y^{\delta_L - (n_1 + m_1 + 2 m_2 +\theta+ 6)}.\] 
It is shown in \cite{BM13} that 
\begin{align*}
\ov{\mc{P}\D}_5 &= \PP \D_5 \cup \ov{\PP \D}_6 \cup \ov{\PP \E}_6.  
\end{align*}
We now define a section of the following bundle 
\begin{align}
\Psi_{\mathcal{P} \E_6}:  \ov{\mc{P}\D}_5 \cap \mu \lra \mathbb{L}_{\mathcal{P} \E_6} & := 
\hat{\gamma}^{*} \otimes (T\X/\hat{\gamma})^{* 2} \otimes \gamma_{\D}^* \otimes L, \qquad 
\textnormal{given by} \nonumber \\ 
\{\Psi_{\mathcal{P} \E_6}([f], l_{q})\}(v\otimes w^{\otimes 2} \otimes f) & := \nabla^3 f|_{q}(v,w,w) \qquad 
\forall ~~v \in l_q, ~~w \in T_q\X/l_q.
\label{psi_pe6_section_defn}
\end{align}
If the line bundle $L$ is sufficiently $4$-ample, then the section $\Psi_{\mathcal{P} \E_6}$ 
vanishes on $\PP \E_6 \cap \mu$ transversally. Moreover, it does not vanish 
on $\PP \D_6 \cap \mu$. 
Hence 
\begin{align*}
\N(\mc{P}\E_6, n_1, m_1, m_2, \theta)  &= \big\lan e(\mathbb{L}_{\mathcal{P} \E_6}), ~~[\ov{\mc{P}\D}_5]\cap [\mu] \big\ran.  
\end{align*}
This proves \eqref{pe6_chern_eqn}. \qed

\begin{prp}
The number $\N(\mc{P}\E_7, n_1, m_1, m_2, \theta)$   is given by 
\begin{align}
\N(\mc{P}\E_7, n_1, m_1, m_2, \theta) &=  4\N(\mc{P}\E_6, n_1, m_1, m_2, \theta+1) 
+\N(\mc{P}\E_6, n_1, m_1, m_2, \theta) \nonumber \\  
& ~~ + \N(\mc{P}\E_6, n_1+1, m_1, m_2, \theta).  \label{pe7_chern_eqn}
\end{align}
\end{prp}
\pf Let $\mu$ be a generic pseudocycle representing the homology class Poincar\'{e} dual 
to 
\[c_1^{n_1} x_1^{m_1} x_2^{m_2} \lambda^{\theta} y^{\delta_L - (n_1 + m_1 + 2 m_2 +\theta+ 7)}.\] 
It is shown in \cite{BM13} that 
\begin{align*}
\ov{\mc{P}\E}_6 &= \PP \E_6 \cup \ov{\PP \E}_7 \cup \ov{\hat{\mathcal{X}}_8^{\#}}. 
\end{align*}
We now define a section of the following bundle 
\begin{align}
\Psi_{\mathcal{P} \E_7}:  \ov{\mc{P}\E}_6 \cap \mu \lra \mathbb{L}_{\mathcal{P} \E_7} & := 
\hat{\gamma}^{* 4} \otimes \gamma_{\D}^* \otimes L, \qquad 
\textnormal{given by} \nonumber \\ 
\{\Psi_{\mathcal{P} \E_7}([f], l_{q})\}(v^{\otimes 4} \otimes f) & := \nabla^4 f|_{q}(v,v,v,v) \qquad 
\forall ~~v \in l_q.
\label{psi_pe7_section_defn}
\end{align}
If the line bundle $L$ is sufficiently $5$-ample, then the section $\Psi_{\mathcal{P} \E_7}$ 
vanishes on $\PP \E_7 \cap \mu$ transversally. Moreover, it does not vanish on 
$\hat{\mathcal{X}}_8^{\#} \cap \mu$.
Hence 
\begin{align*}
\N(\mc{P}\E_7, n_1, m_1, m_2, \theta)  &= \big\lan e(\mathbb{L}_{\mathcal{P} \E_7}), ~~[\ov{\mc{P}\E}_6] \cap [\mu] \big\ran.  
\end{align*}
This proves \eqref{pe6_chern_eqn}. \qed 

\begin{prp}
The number $\N(\mc{P}\D_6, n_1, m_1, m_2, \theta)$   is given by 
\begin{align}
\N(\mc{P}\D_6, n_1, m_1, m_2, \theta) &=  4\N(\mc{P}\D_5, n_1, m_1, m_2, \theta+1) 
+\N(\mc{P}\D_5, n_1, m_1+1, m_2, \theta) \nonumber \\  
& ~~ + \N(\mc{P}\D_5, n_1+1, m_1, m_2, \theta).  \label{pd6_chern_eqn}
\end{align}
\end{prp}
\pf Let $\mu$ be a generic pseudocycle representing the homology class Poincar\'{e} dual 
to 
\[c_1^{n_1} x_1^{m_1} x_2^{m_2} \lambda^{\theta} y^{\delta_L - (n_1 + m_1 + 2 m_2 + \theta+6)}.\] 
It is shown in \cite{BM13} that 
\begin{align*}
\ov{\mc{P}\D}_5 &= \PP \D_5 \cup \ov{\PP \D}_6 \cup \ov{\PP \E}_6.  
\end{align*}
We now define a section of the following bundle 
\begin{align}
\Psi_{\mathcal{P} \D_6}:  \ov{\mc{P}\D}_5 \cap \mu \lra \mathbb{L}_{\mathcal{P} \D_6} & := 
\hat{\gamma}^{*4} \otimes \gamma_{\D}^* \otimes L, \qquad 
\textnormal{given by} \nonumber \\ 
\{\Psi_{\mathcal{P} \D_6}([f], l_{q})\}(v^{\otimes 4}\otimes f) & := \nabla^4 f|_{q}(v,v,v,v) \qquad 
\forall ~~v \in \hat{\gamma}. \label{psi_pd6_section_defn}
\end{align}
If the line bundle $L$ is sufficiently $5$-ample, then the section $\Psi_{\mathcal{P} \D_6}$ 
vanishes on $\PP \D_6 \cap \mu$ transversally. Moreover, it does not vanish on $\PP \E_6 \cap \mu$. 
Hence 
\begin{align*}
\N(\mc{P}\D_6, n_1, m_1, m_2, \theta)  &= \big\lan e(\mathbb{L}_{\mathcal{P} \D_6}), ~~[\ov{\mc{P}\D}_5] \cap [\mu] \big\ran.  
\end{align*}
This proves \eqref{pd6_chern_eqn}. \qed \\[0.2cm]
\ni Let us now make a small abbreviation: 
if $v \in l_q$ and $w \in TX/l_q$, we define  
\begin{align*}
f_{ij} & := \nabla^{i+j} f|_{q}
(\underbrace{v,\cdots v}_{\textnormal{$i$ times}}, \underbrace{w,\cdots w}_{\textnormal{$j$ times}}).
\end{align*}

\begin{prp}
The number $\N(\mc{P}\D_7, n_1, m_1, m_2, \theta)$   is given by 
\begin{align}
\N(\mc{P}\D_7, n_1, m_1, m_2, \theta) &=  4\N(\mc{P}\D_6, n_1, m_1, m_2, \theta+1) 
+2\N(\mc{P}\D_6, n_1, m_1+1, m_2, \theta) \nonumber \\  
& ~~ + 2\N(\mc{P}\D_6, n_1, m_1, m_2, \theta) + 2\N(\mc{P}\D_6, n_1+1, m_1, m_2, \theta).  \label{pd7_chern_eqn}
\end{align}
\end{prp}
\pf Let $\mu$ be a generic pseudocycle representing the homology class Poincar\'{e} dual 
to 
\[c_1^{n_1} x_1^{m_1} x_2^{m_2} \lambda^{\theta} y^{\delta_L - (n_1 + m_1 + 2 m_2 +\theta+ 7)}.\] 
It is shown in \cite{BM13} that 
\begin{align}
\ov{\mc{P}\D}_6 &= \mc{P}\D_6 \cup \ov{\mc{P}\D}_7 \cup \ov{\mc{P}\E}_7.   
\end{align}
We now define a section of the following bundle 
\begin{align}
\Psi_{\mathcal{P} \D_7}:  \ov{\mc{P}\D}_6 \cap \mu \lra \mathbb{L}_{\mathcal{P} \D_7} & := 
\hat{\gamma}^{*6} \otimes (TX/\hat{\gamma})^{*2} \otimes \gamma_{\D}^{*2} \otimes L^{\otimes 2}, \qquad 
\textnormal{given by} \nonumber \\ 
\{\Psi_{\mathcal{P} \D_7}([f], l_{q})\}(v^{\otimes 6} \otimes w^{\otimes 2} \otimes f^{\otimes 2}) 
& := f_{12}\D^{f}_7 \qquad 
\forall ~~v \in \hat{\gamma},  ~~w \in TX/\hat{\gamma}, 
\end{align}
where $~\D^{f}_7 := f_{50} - \frac{5 f_{31}^2}{3 f_{12}}.$
If the line bundle $L$ is sufficiently $6$-ample, then the section $\Psi_{\mathcal{P} \D_7}$, 
restricted to $\mc{P} \D_7 \cap \mu$ vanishes transversally.  
Moreover, the section does not vanish on $\mc{P}\E_7 \cap \mu$. 
Hence 
\begin{align*}
\N(\mc{P}\D_7, n_1, m_1, m_2, \theta)  &= \big\lan e(\mathbb{L}_{\mathcal{P} \D_7}), ~~[\ov{\mc{P}\D}_6]\cap [\mu] \big\ran.  
\end{align*}
This proves \eqref{pd7_chern_eqn}. \qed \\

\begin{prp}
The number $\N(\mc{P}\A_4, n_1, m_1, m_2, \theta)$   is given by 
\begin{align}
\N(\mc{P}\A_4, n_1, m_1, m_2, \theta) &=  2\N(\mc{P}\A_3, n_1, m_1, m_2, \theta+1) 
+2\N(\mc{P}\A_3, n_1, m_1+1, m_2, \theta) \nonumber \\  
& ~~ + 2 \N(\mc{P}\A_3, n_1, m_1, m_2, \theta) + 2\N(\mc{P}\A_3, n_1+1, m_1, m_2, \theta).  \label{pa4_chern_eqn}
\end{align}
\end{prp}
\pf Let $\mu$ be a generic pseudocycle representing the homology class Poincar\'{e} dual 
to 
\[c_1^{n_1} x_1^{m_1} x_2^{m_2} \lambda^{\theta} y^{\delta_L - (n_1 + m_1 + 2 m_2 +\theta+ 4)}.\] 
It is shown in \cite{BM13} that 
\begin{align}
\ov{\mc{P}\A}_3 &= \mc{P}\A_3 \cup \ov{\mc{P}\A}_4 \cup \ov{\mc{P}\D}_4.   
\end{align}
We now define a section of the following bundle 
\begin{align}
\Psi_{\mathcal{P} \A_4}:  \ov{\mc{P}\A}_3 \cap \mu \lra \mathbb{L}_{\mathcal{P} \A_4} & := 
\hat{\gamma}^{*4} \otimes (TX/\hat{\gamma})^{*2} \otimes \gamma_{\D}^{*2} \otimes L^{\otimes 2}, \qquad 
\textnormal{given by} \nonumber \\ 
\{\Psi_{\mathcal{P} \A_4}([f], l_{q})\}(v^{\otimes 4} \otimes w^{\otimes 2} \otimes f^{\otimes 2}) 
& := f_{02}\A^{f}_4 \qquad 
\forall ~~v \in l_q,  ~~w \in T_qX/l_q, 
\end{align}
where $~\A^{f}_4 := f_{40} - \frac{3 f_{21}^2}{f_{02}}.$
If the line bundle $L$ is sufficiently $5$-ample, then the section $\Psi_{\mathcal{P} \A_4}$, 
restricted to $\mc{P} \A_4 \cap \mu$ vanishes transversally.  
Moreover, the section does not vanish on $\mc{P}\D_4 \cap \mu$.
Hence 
\begin{align*}
\N(\mc{P}\A_4, n_1, m_1, m_2, \theta)  &= \big\lan e(\mathbb{L}_{\mathcal{P} \A_4}), ~~[\ov{\mc{P}\A}_3] \big\ran.  
\end{align*}
This proves \eqref{pa4_chern_eqn}. \qed \\

\begin{prp}
The number $\N(\mc{P}\A_5, n_1, m_1, m_2, \theta)$   is given by 
\begin{align}
\N(\mc{P}\A_5, n_1, m_1, m_2, \theta) &=  \N(\mc{P}\A_4, n_1, m_1, m_2, \theta+1) 
+4\N(\mc{P}\A_4, n_1, m_1+1, m_2, \theta) \nonumber \\  
& ~~ + 3\N(\mc{P}\A_4, n_1, m_1, m_2, \theta) + 3\N(\mc{P}\A_4, n_1+1, m_1, m_2, \theta) \nonumber  \\
& ~~ -2\N(\mc{P}\D_5,n_1,m_1,m_2, \theta).  \label{pa5_chern_eqn}
\end{align}
\end{prp}
\pf Let $\mu$ be a generic pseudocycle representing the homology class Poincar\'{e} dual 
to 
\[c_1^{n_1} x_1^{m_1} x_2^{m_2} \lambda^{\theta} y^{\delta_L - (n_1 + m_1 + 2 m_2 +\theta+ 5)}.\] 
It is shown in \cite{BM13} that 
\begin{align}
\ov{\mc{P}\A}_4 &= \mc{P}\A_4 \cup \ov{\mc{P}\A}_5 \cup \ov{\mc{P}\D}_5.   
\end{align}
We now define a section of the following bundle 
\begin{align}
\Psi_{\mathcal{P} \A_5}:  \ov{\mc{P}\A}_4 \cap \mu \lra \mathbb{L}_{\mathcal{P} \A_5} & := 
\hat{\gamma}^{*5} \otimes (TX/\hat{\gamma})^{*4} \otimes \gamma_{\D}^{*3} \otimes L^{\otimes 3}, \qquad 
\textnormal{given by} \nonumber \\ 
\{\Psi_{\mathcal{P} \A_4}([f], l_{q})\}(v^{\otimes 5} \otimes w^{\otimes 4} \otimes f^{\otimes 3}) 
& := f_{02}^2\A^{f}_5 \qquad 
\forall ~~v \in l_q,  ~~w \in T_qX/l_q, 
\end{align}
where 
\[ \A^{f}_5 := f_{50} -\frac{10 f_{21} f_{31}}{f_{02}} + 
\frac{15 f_{12} f_{21}^2}{f_{02}^2}.\]
If the line bundle $L$ is sufficiently $6$-ample, then the section $\Psi_{\mathcal{P} \A_5}$, 
restricted to $\mc{P} \A_5$ vanishes transversally. It is also shown in   
\cite{BM13}, that this section vanishes on $\mc{P}\D_5 \cap \mu$ with a multiplicity of $2$.  
Hence 
\begin{align*}
\big\lan e(\mathbb{L}_{\mathcal{P} \A_5}), ~~[\ov{\mc{P}\A}_4] \cap [\mu] \big\ran &=
\N(\mc{P}\A_5, n_1, m_1, m_2, \theta) + 
2\N(\mc{P}\D_5, n_1, m_1, m_2, \theta).     
\end{align*}
This proves \eqref{pa4_chern_eqn}. \qed \\

\begin{prp}
The number $\N(\mc{P}\A_6, n_1, m_1, m_2, \theta)$   is given by 
\begin{align}
\N(\mc{P}\A_6, n_1, m_1, m_2, \theta) &=  6\N(\mc{P}\A_5, n_1, m_1+1, m_2, \theta) \nonumber \\  
& ~~ + 4\N(\mc{P}\A_5, n_1, m_1, m_2, \theta) + 4\N(\mc{P}\A_5, n_1+1, m_1, m_2, \theta) \nonumber  \\
& ~~ -4\N(\mc{P}\D_6,n_1,m_1,m_2, \theta)-3\N(\mc{P}\E_6,n_1,m_1,m_2, \theta).  \label{pa6_chern_eqn}
\end{align}
\end{prp}
\pf Let $\mu$ be a generic pseudocycle representing the homology class Poincar\'{e} dual 
to 
\[c_1^{n_1} x_1^{m_1} x_2^{m_2} \lambda^{\theta} y^{\delta_L - (n_1 + m_1 + 2 m_2 +\theta+ 6)}.\] 
It is shown in \cite{BM13} that 
\begin{align}
\ov{\mc{P}\A}_5 &= \mc{P}\A_5 \cup \ov{\mc{P}\A}_6 \cup \ov{\mc{P}\D}_6 \cup \ov{\mc{P}\E}_6.   
\end{align}
We now define a section of the following bundle 
\begin{align}
\Psi_{\mathcal{P} \A_6}:  \ov{\mc{P}\A}_5 \lra \mathbb{L}_{\mathcal{P} \A_6} & := 
\hat{\gamma}^{*6} \otimes (TX/\hat{\gamma})^{*6} \otimes \gamma_{\D}^{*4} \otimes L^{\otimes 4}, \qquad 
\textnormal{given by} \nonumber \\ 
\{\Psi_{\mathcal{P} \A_4}([f], l_{q})\}(v^{\otimes 6} \otimes w^{\otimes 6} \otimes f^{\otimes 4}) 
& := f_{02}^3\A^{f}_6 \qquad 
\forall ~~v \in l_q,  ~~w \in T_q\X/l_q, 
\end{align}
where 
\[ \A^{f}_6 := f_{60}- \f{ 15 f_{21} f_{41}}{f_{02}}-\f{10 f_{31}^2}{f_{02}} + \f{60 f_{12} f_{21} f_{31}}{f_{02}^2}
   +
   \f{45 f_{21}^2 f_{22}}{f_{02}^2} - \f{15 f_{03} f_{21}^3}{f_{02}^3}
   -\f{90 f_{12}^2 f_{21}^2}{f_{02}^3} .\]
If the line bundle $L$ is sufficiently $7$-ample, then the section $\Psi_{\mathcal{P} \A_6}$, 
restricted to $\mc{P} \A_6$ vanishes transversally. It is shown in   
\cite{BM13}, that this section  vanishes on $\mc{P}\D_6 \cap \mu$ and 
$\mc{P}\E_6 \cap \mu$
with a multiplicity of $4$ and $3$ respectively.  
Hence 
\begin{align*}
\lan e(\mathbb{L}_{\mathcal{P} \A_6}), [\ov{\mc{P}\A}_5] \cap [\mu] \ran &=
\N(\mc{P}\A_6, n_1, m_1, m_2, \theta) + 
4\N(\mc{P}\D_6, n_1, m_1, m_2, \theta)+3\N(\mc{P}\E_6, n_1, m_1, m_2, \theta).     
\end{align*}
This proves \eqref{pa4_chern_eqn}. \qed \\

\begin{prp}
\label{pa7_chern}
The number $\N(\mc{P}\A_7, n_1, m_1, m_2, 0)$   is given by 
\begin{align}
\N(\mc{P}\A_7, n_1, m_1, m_2, 0) &=  -\N(\mc{P}\A_6, n_1, m_1, m_2, 1) 
+ 8 \N(\mc{P}\A_6, n_1, m_1+1, m_2, 0) \nonumber \\  
& ~~ + 5\N(\mc{P}\A_6, n_1, m_1, m_2, 0) + 5\N(\mc{P}\A_5, n_1+1, m_1, m_2, 0) \nonumber  \\
& ~~ -6\N(\mc{P}\D_6,n_1,m_1,m_2, 0)-7\N(\mc{P}\E_6,n_1,m_1,m_2, 0).  \label{pa7_chern_eqn}
\end{align}
\end{prp}
\pf Let $\mu$ be a generic pseudocycle representing the homology class Poincar\'{e} dual 
to 
\[c_1^{n_1} x_1^{m_1} x_2^{m_2} y^{\delta_L - (n_1 + m_1 + 2 m_2 + 7)}.\] 
It is shown in \cite{BM13} that 
\begin{align}
\ov{\mc{P}\A}_6 &= \mc{P}\A_6 \cup \ov{\mc{P}\A}_7 \cup \ov{\mc{P}\D}_7 \cup \ov{\mc{P}\E}_7 \cup 
\ov{\hat{\mathcal{X}}}_8.   
\end{align}
We now define a section of the following bundle 
\begin{align}
\Psi_{\mathcal{P} \A_7}:  \ov{\mc{P}\A}_6 \cap \mu \lra \mathbb{L}_{\mathcal{P} \A_7} & := 
\hat{\gamma}^{*7} \otimes (TX/\hat{\gamma})^{*8} \otimes \gamma_{\D}^{*5} \otimes L^{\otimes 5}, \qquad 
\textnormal{given by} \nonumber \\ 
\{\Psi_{\mathcal{P} \A_7}([f], l_{q})\}(v^{\otimes 7} \otimes w^{\otimes 8} \otimes f^{\otimes 5}) 
& := f_{02}^4\A^{f}_7 \qquad 
\forall ~~v \in l_q,  ~~w \in T_q\X/l_q, 
\end{align}
where 
\begin{align*}
\A^{f}_7 & := f_{70} - \frac{21 f_{21} f_{51}}{f_{02}} 
- \frac{35 f_{31} f_{41}}{f_{02}} + \frac{105 f_{12} f_{21} f_{41}}{f_{02}^2} + \f{105 f_{21}^2 f_{32}}{f_{02}^2} + 
\f{70 f_{12} f_{31}^2}{f_{02}^2}+ \f{210 f_{21}f_{22}f_{31}}{f_{02}^2} \nonumber \\
&
-\f{105 f_{03} f_{21}^2 f_{31}}{ f_{02}^3}
-\f{420 f_{12}^2 f_{21} f_{31}}{f_{02}^3}
-\f{630 f_{12}f_{21}^2 f_{22}}{f_{02}^3}
-\f{105 f_{13} f_{21}^3}{f_{02}^3}
+ \f{315 f_{03} f_{12} f_{21}^3}{f_{02}^4}
+ \f{630 f_{12}^3 f_{21}^2}{f_{02}^4}.
\end{align*}
If the line bundle $L$ is sufficiently $8$-ample, then the section $\Psi_{\mathcal{P} \A_7}$, 
restricted to $\mc{P} \A_7 \cap \mu$ vanishes transversally. It is shown in   
\cite{BM13}, that this section  vanishes on $\mc{P}\D_7 \cap \mu$ and 
$\mc{P}\E_7 \cap \mu$
with a multiplicity of $6$ and $7$ respectively. Assume that it vanishes on $\hat{\mc{X}}_8$ with a 
multiplicity of $\eta$. 
Hence 
\begin{align*}
\lan e(\mathbb{L}_{\mathcal{P} \A_7}), ~~[\ov{\mc{P}\A}_6] \cap [\mu] \ran &=
\N(\mc{P}\A_7, n_1, m_1, m_2, 0) + 
6\N(\mc{P}\D_7, n_1, m_1, m_2, 0)+7\N(\mc{P}\E_7, n_1, m_1, m_2, 0) \\ 
& + \eta |\ov{\hat{\mc{X}}}_8 \cap \mu|.
\end{align*}
Since $|\ov{\hat{\mc{X}}}_8 \cap \mu| =0$, we get \eqref{pa7_chern_eqn}. \qed \\

\section{Topological computations: two singular points} 
\label{top_comp_two_sing_pt}
\hf\hf We will now give a proof of Theorem \ref{two_pt_chern_class_codim_7}.
Before that, let us setup some additional notation. 
Let $S_1$ and $S_2$ be two subsets of $\D \times \X$ and $T_2$ be a subset of $\D \times \mathbb{P}T\X$.
Define 
\begin{align}
S_1 \circ S_2 &:= 
\{ ([f], q_1, q_2) \in \D \times \X \times \X:  ([f], q_1) \in S_1, ~~([f], q_2) \in S_2, ~~q_1 \neq q_2 \},
\nonumber \\
S_1 \circ T_2 &:= 
\{ ([f], q_1, l_{q_2}) \in \D \times \X \times \mathbb{P} T\X:  ([f], q_1) \in S_1, ~~([f], l_{q_2}) \in T_2, 
~~q_1 \neq q_2 \}. \nonumber 
\end{align}
Next, if $\mf{X}$ is a codimension $k$ singularity, 
we define the following numbers: 
\begin{align*}
\N(\A_1\mf{X}, n_1, m_1, m_2) &:= 
\lan c_1^{n_1} x_1^{m_1} x_2^{m_2} y^{\delta_L-(n_1+m_1+2m_2 + k+1)},  
~[\overline{\A_1 \circ \mf{X}}]\ran,  \\
\N(\A_1 \mc{P} \mf{X}, n_1, m_1, m_2, \theta) &:= \lan c_1^{n_1} x_1^{m_1} x_2^{m_2} 
\lambda^{\theta} y^{\delta_L-(n_1+m_1+2m_2 +\theta+ k+1)},  
~[\overline{\A_1 \circ \mc{P} \mf{X}}]\ran.   
\end{align*}
Next, if $S$ and $T$ are subsets of  $\D \times \X$ and 
$\D \times \mathbb{P}T\X$ respectively, we define 
\begin{align*}
\Delta S &:= \{([f], q, q) \in \D \times \X \times \X:  ([f], q) \in S \} \qquad \textnormal{and} \\
\Delta T & := \{([f], q, l_q) \in \D \times \X \times \mathbb{P}T\X:  ([f], l_q) \in T \}. 
\end{align*}
Finally, 
we    
define the following two projection maps
\[ \pi_2:=\textup{id}\times\textup{proj}_2: \D \times \X \times \X \lra \D \times \X, \qquad \qquad \pi_2:=\textup{id}\times\textup{proj}_2: \D \times \X \times \mathbb{P} T\X \lra \D \times \mathbb{P} T\X \]
where $\textup{proj}_2$ denotes the projection onto the second factor. \\
\hf \hf Note that given any bundle over $\D \times \X$ (resp. $\D \times \P T\X$), there is 
an induced bundle over $\D \times \X \times \X$ (resp. $\D \times \X \times \P T\X$) arising from the 
pullback via $\pi_2$. Similarly, given any section of such a bundle, there is a corresponding 
section on the pullback bundle, via $\pi_2$. We will encounter 
the bundles and sections of these bundles that 
we defined in section \ref{top_comp_one_sing_pt}; we will encounter them over 
$\D \times \X \times \X$ or $\D \times \X \times \P T\X$ via the pullback of $\pi_2$. \\
\hf \hf  We will now give a series of formulas to compute $\N(\A_1\A_1, n_1, m_1, m_2)$ and 
$\N(\A_1\mc{P}\mf{X}, n_1, m_1, m_2, \theta)$. 
Note that $\N(\A_1\A_1) = \N(\A_1\A_1, 0,0,0)$.
In order to compute the remaining $\N(\A_1\mf{X})$ we 
do the following: if $\mf{X}$ is anything other than $\D_4$, then we observe that 
$\N(\A_1\mf{X}) = \N(\A_1\mc{P}\mf{X}, 0, 0, 0, 0)$. If $\mf{X} = \D_4$, then we observe that 
$\N(\A_1\D_4) = \frac{\N(\A_1\mc{P}\D_4,0,0,0,0)}{3}$.\\  
\hf \hf We are now ready to give a proof of Theorem \ref{two_pt_chern_class_codim_7}. 
An important notational remark is that in the subsequent proofs we shall use $\mu$ to 
denote a homology class Poincar\'{e} dual to an 
element $c_1^{n_1} x_1^{m_1} x_2^{m_2}\lambda^\theta y^k$ in $H^*(\D\times X\times \P T\X ;\mathbb{Z})$. 
The elements $\lambda^\theta$ and $y^k$ are unambiguously defined via the appropriate pullbacks. 
But the class $c_1^{n_1}x_1^{m_1}x_2^{m_2}$ is understood to arise from the second 
factor in $X\times \mathbb{P} T\X$, i.e., if $\pi_2:\D\times X\times \P T\X \to \D \times \P T\X$ 
is the projection map defined earlier, then 
\bgd
c_1:=c_1(\pi_2^\ast(L)),\,\,x_i=c_i(\pi_2^*(T^*X)).
\edd
We will allow ourselves this abuse of notation. 

\begin{prp}
\label{a1a1_chern}
The number $\N(\A_1\A_1, n_1, m_1, m_2)$ is given by 
\begin{align}
\label{a1a1_chern_eqn}
\N(\A_1\A_1, n_1, m_1, m_2) & = \N(\A_1) \times \N(\A_1,n_1, m_1, m_2) \nonumber \\ 
                            &   -\Big(\N(\A_1,n_1, m_1, m_2) + \N(\A_1, n_1+1, m_1, m_2) 
                            + 3 \N(\A_2, n_1, m_1, m_2) \Big).
\end{align}
\end{prp}
\pf Let $\mu$ be a pseudocycle 
in $\D \times \X \times \X$
representing the homology class 
Poincar\'{e} dual to  
\[c_1^{n_1} x_1^{m_1} x_2^{m_2} y^{\delta_L-(n_1+m_1+2m_2 + 2)}.\]
Note that 
\begin{align*} 
\ov{\A}_1 \times \X & = \ov{\ov{\A}_1 \circ (\D \times \X)} = \ov{\A}_1 \circ (\D \times \X) \du \Delta \ov{\A}_1.
\end{align*} 
We show in \cite{BM13_2pt_published} and \cite{BM_Detail} that 
the sections 
\begin{align}
 \pi_2^*\ds_{\A_0}:\ov{\A}_1 \times \X - \Delta \ov{\A}_1   \lra \pi_2^* \DL_{\A_0}, 
 \qquad \pi_2^*\ds_{\A_1}: \pi_2^*\ds_{\A_0}^{-1}(0) \lra \pi_2^*\DV_{\A_1}
\end{align}
are transverse to the zero set (if $L$ is sufficiently $4$-ample). 
Hence 
\begin{align}
\big\langle e(\pi_2^*\DL_{\A_0}) e(\pi_2^*\DV_{\A_1}), ~[\ov{\A}_1 \times \X] \cap [\mu]  \big\rangle & =  
\Num(\A_1 \A_1, n_1, m_1, m_2) + 
\mathcal{C}_{\Delta \ov{\A}_1 \cap \mu}(\pi_2^*\ds_{\A_0} \oplus \pi_2^*\ds_{\A_1}), \label{A1A1_Euler_Class}
\end{align}
where $\mathcal{C}_{\Delta \ov{\A}_1 \cap \mu}(\pi_2^*\ds_{\A_0} \oplus \pi_2^*\ds_{\A_1})$ 
is the contribution of the section 
$\pi_2^*\ds_{\A_0} \oplus \pi_2^*\ds_{\A_1}$ to the Euler class from 
$\Delta \ov{\A}_1 \cap \mu$. 
The lhs of \eqref{A1A1_Euler_Class}, as computed by splitting principle and a case by case check, is 
\begin{align}
\big\langle e(\pi_2^*\DL_{\A_0}) e(\pi_2^*\DV_{\A_1}), ~[\ov{\A}_1 \times \X]\cap[\mu] \big\rangle  &= 
\N(\A_1) \times \N(\A_1,n_1, m_1, m_2) \label{a1a1_Euler_class_Main_stratum}.
\end{align}
\ni Next, we  compute 
$\mathcal{C}_{\Delta \ov{\A}_1 \cap \mu}(\pi_2^*\ds_{\A_0}\oplus\pi_2^*\ds_{\A_1})$.
Note that ~$ \ov{\A}_1 = \A_1 \du \ov{\A}_2.$
It is shown in \cite{BM13_2pt_published} that 
\begin{align}
\mathcal{C}_{\Delta \A_1 \cap \mu}(\pi_2^*\ds_{\A_0}\oplus\pi_2^*\ds_{\A_1}) & = 
\big\langle e(\pi_2^* \DL_{\A_0}) , ~[\Delta \ov{\A}_1] \cap [\mu] \big\rangle \nonumber \\ 
& = \N(\A_1,n_1, m_1, m_2) + \N(\A_1, n_1+1, m_1, m_2),  \label{a1a1_a1_contribution}\\
\mathcal{C}_{\Delta \ov{\A}_2 \cap \mu}(\pi_2^*\ds_{\A_0}\oplus\pi_2^*\ds_{\A_1}) &= 
3 \Num(\A_2, n_1, m_1, m_2).  \label{a1a1_a2_contribution}
\end{align}
It is easy to see that \eqref{a1a1_Euler_class_Main_stratum}, \eqref{a1a1_a1_contribution} and \eqref{a1a1_a2_contribution} 
prove \eqref{a1a1_chern_eqn}. \qed \\

\begin{prp}
\label{a1pa2_chern}
The number $\N(\A_1\mc{P}\A_2, n_1, m_1, m_2, \theta)$ is given by 
\begin{align}
\label{a1pa2_chern_eqn}
\N(\A_1\mc{P}\A_2, n_1, m_1, m_2, \theta) &= 
\begin{cases} 2\N(\A_1\A_1, n_1, m_1, m_2) + 2 \N(\A_1\A_1, n_1, m_1+1, m_2) \\ 
+ 2\N(\A_1\A_1, n_1+1, m_1, m_2)  
 -\Big(2 \N(\mc{P}\A_3, n_1, m_1, m_2, \theta) \Big) \qquad 
\mbox{if} ~~\theta =0,\\ 
\\
                                  \N(\A_1\A_1, n_1, m_1, m_2) + 2\N(\A_1\A_1, n_1+1, m_1, m_2) + 
                                  \N(\A_1\A_1, n_1+2, m_1, m_2) \\ 
                                  + 3 \N(\A_1\A_1, n_1, m_1+1, m_2) + 3 \N(\A_1\A_1, n_1+1, m_1+1, m_2) \\ 
                                  + 
                                  2\N(\A_1\A_1, n_1, m_1+2, m_2)\\  
                                  -\Big( 2 \N(\mc{P}\A_3, n_1, m_1, m_2, \theta) + 3 \N(\mc{D}_4, n_1, m_1, m_2) \Big)
                                  \qquad \mbox{if} ~~\theta =1,\\ \\
                                  \N(\A_1\mc{P}\A_2, n_1, m_1+1, m_2, \theta-1) \\
                                 -\N(\A_1\mc{P}\A_2, n_1, m_1, m_2+1, \theta-2) 
                                 \qquad \mbox{if} ~~\theta >1. \end{cases}
\end{align}
\end{prp}
\pf Let $\mu$ be a pseudocycle 
in $\D \times \X \times \P T\X$
representing the homology class 
Poincar\'{e} dual to 
\[ c_1^{n_1} x_1^{m_1} x_2^{m_2} \lambda^{\theta} y^{\delta_L-(n_1+m_1+2m_2 + \theta+3)}.\] 
We show in \cite{BM13_2pt_published}, that 
\begin{align*}
\ov{\ov{\A}_1 \circ  \hat{\A}^{\#}_1}& = 
\ov{\A}_1\circ \hat{\A}^{\#}_1  \du \ov{\A}_1 
                           \circ  \ov{\PP \A}_2 \du \Delta \ov{\hat{\A}}_3.
\end{align*}
If $L$ is sufficiently $5$-ample, then 
the section 
\[ \pi_2^*\us_{\PP \A_2} : \ov{\ov{\A}_1 \circ  \hat{\A}^{\#}_1} \cap \mu \lra \pi_2^*\UV_{\PP \A_2}\]
vanishes on $\A_1 \circ \PP \A_2 \cap \mu$ transversely. 
Hence, the zeros of the section 
$$ \pi_2^*\us_{\PP \A_2} : \ov{\ov{\A}_1 \circ  \hat{\A}^{\#}_1} \cap \mu 
\lra \pi_2^*\UV_{\PP \A_2}, $$
restricted to $\A_1 \circ \PP \A_2 \cap \mu$ 
counted with a sign, is our desired number. 
In other words 
\bgd
\Big\langle e(\pi_2^*\UV_{\PP \A_2}), ~[\ov{\ov{\A}_1 \circ  \hat{\A}^{\#}_1}] \cap \mu \Big\rangle = 
\Num(\A_1 \PP\A_2, n_1,m_1,m_2,\theta) + \mathcal{C}_{\Delta \ov{\hat{\A}}_3 \cap \mu}\Big(\pi_2^*\us_{\PP \A_2}\Big) 
\edd
where $\mathcal{C}_{\Delta \ov{\hat{\A}}_3 \cap \mu}\Big(\pi_2^*\us_{\PP \A_2}\Big)$ is the 
contribution of the section to the Euler class from 
$\Delta \ov{\hat{\A}}_3 \cap \mu$. 
Note that $ \pi_2^*\us_{\PP \A_2}$ vanishes only on $\Delta \PP \A_3 \cap \mu$ and $\Delta \hat{\D}_4^{\#\flat} \cap \mu$ and 
not on the entire $\Delta \ov{\hat{\A}}_3 \cap \mu$.   
We show in \cite{BM13_2pt_published} that 
the contribution from $\Delta \PP \A_3 \cap \mu$ and $\Delta \hat{\D}_4^{\#\flat} \cap \mu$ are $2$ and $3$ respectively. 
Hence 
\bgd
\Big\lan e(\pi_2^*\UV_{\PP \A_2}), ~~[\ov{\ov{\A}_1 \circ  \hat{\A}^{\#}_1}] \cap \mu \Big\ran = 
\Num(\A_1\PP\A_2,n_1,m_1, m_2, \theta) + 2 \Num(\PP \A_3, n_1, m_1, m_2, \theta) + 3 |\Delta \ov{\hat{\D}_4^{\#\flat}} \cap \mu|.  
\edd
Since $\ov{\hat{\D}_4^{\#\flat}}=\ov{\hat{\D}_4}$, this gives us \eqref{a1pa2_chern_eqn}.   \qed \\

\begin{prp}
\label{a1pa3_chern}
The number $\N(\A_1\mc{P}\A_3, n_1, m_1, m_2, \theta)$ is given by 
\begin{align}
\N(\A_1\mc{P}\A_3, n_1, m_1, m_2, \theta) &= 3\N(\A_1\mc{P}\A_2, n_1, m_1, m_2, \theta+1) 
+\N(\A_1\mc{P}\A_2, n_1, m_1, m_2, \theta) \nonumber \\  
& ~~ + \N(\A_1\mc{P}\A_2, n_1+1, m_1, m_2, \theta) -2 \N(\mc{P}\A_5, n_1, m_1, m_2, \theta).  \label{a1pa3_chern_eqn}
\end{align}
\end{prp}
\pf Let $\mu$ be a generic pseudocycle 
representing the homology class 
Poincar\'{e} dual to 
\[ c_1^{n_1} x_1^{m_1} x_2^{m_2} \lambda^{\theta} y^{\delta_L-(n_1+m_1+2m_2 + \theta+4)}.\] 
We show in \cite{BM13_2pt_published}, that 
\begin{align*}
\ov{\ov{\A}_1 \circ  \PP \A}_2
                         &= \ov{\A}_1 \circ  \PP \A_2 \du \ov{\A}_1 \circ (\ov{\PP \A}_3 \cup \ov{\hat{\D}^{\#}_4}) \du 
\Big( \Delta \ov{\PP \A}_4 \cup \Delta \ov{\hat{\D}^{\#\flat}_5}\Big).
\end{align*}
If $L$ is sufficiently $6$-ample, then 
the section 
\[ \pi_2^*\us_{\PP \A_3} : \ov{\ov{\A}_1 \circ \PP \A}_2 \cap \mu \lra \pi_2^*\UL_{\PP \A_3}\]
vanishes transversely on $\A_1 \circ \PP \A_3$. 
We also show in \cite{BM13_2pt_published}, that
the contribution to the Euler class from the points of $\Delta \PP \A_4 \cap \mu$ is $2$. 
This section does not vanish on $\A_1\circ \hat{\D}_4^{\#}$ and by definition it also does not vanish on $\Delta\hat{\D}^{\#\flat}_5$.
Hence 
\bgd
\Big\lan e(\pi_2^*\UL_{\PP \A_3}), ~~[\ov{\ov{\A}_1 \circ \PP \A}_2] \cap \mu \Big\ran = 
\Num(\A_1\PP\A_3,n_1,m_1, m_2) + 2 \Num(\PP \A_4, n_1, m_1, m_2, \theta)  
\edd
which gives us \eqref{a1pa3_chern_eqn}.  \qed \\


\begin{prp}
\label{a1pa4_chern}
The number $\N(\A_1\mc{P}\A_4, n_1, m_1, m_2, \theta)$ is given by 
\begin{align}
\N(\A_1\mc{P}\A_4, n_1, m_1, m_2, \theta) &= 2\N(\A_1\mc{P}\A_3, n_1, m_1, m_2, \theta+1) 
+2 \N(\A_1\mc{P}\A_3, n_1, m_1+1, m_2, \theta) \nonumber \\  
& ~~ + 2 \N(\A_1\mc{P}\A_3, n_1, m_1, m_2, \theta)+2 \N(\A_1\mc{P}\A_3, n_1+1, m_1, m_2, \theta) \nonumber  \\  
& ~~ -2 \N(\mc{P}\A_5, n_1, m_1, m_2, \theta).  \label{a1pa4_chern_eqn}
\end{align}
\end{prp}
\pf Let $\mu$ be a generic pseudocycle 
representing the homology class 
Poincar\'{e} dual to 
\[ c_1^{n_1} x_1^{m_1} x_2^{m_2} \lambda^{\theta} y^{\delta_L-(n_1+m_1+2m_2 + \theta+5)}.\]
We show in \cite{BM13_2pt_published}, that
\begin{align*}
\ov{\ov{\A}_1 \circ  \PP \A}_3 &= 
\ov{\A}_1 \circ  \PP \A_3 \du \ov{\A}_1 \circ (\ov{\PP \A}_4 \cup \ov{\PP \D}_4) \du 
\Big( \Delta \ov{\PP \A}_5 \cup \Delta \ov{\PP \D^{\vee}_5} \Big). 
\end{align*}
If $L$ is sufficiently $7$-ample, then 
the section 
$$ \pi_2^*\us_{\PP \A_4} : \ov{\ov{\A}_1 \circ \PP \A}_3 \cap \mu \lra \pi_2^*\UL_{\PP \A_4} $$
vanishes transversely on $\A_1 \circ \PP \A_4$. 
It is easy to see that it does not vanish on $\A_1 \circ \PP \D_4 \cap \mu$ and 
$\Delta \PP \D^{\vee}_5 \cap \mu$. 
We also show in \cite{BM13_2pt_published} that 
the contribution of this section to the Euler class from the points of $\Delta\PP \A_5 \cap \mu$ is $2$. 
Hence 
\bgd
\Big\langle e(\pi_2^*\UL_{\PP \A_4}), ~~[\ov{\ov{\A}_1 \circ \PP \A}_3] \cap [\mu] \Big\rangle = 
\Num(\A_1\PP\A_4,n_1,m_1,m_2, \theta) + 2 \Num(\PP \A_5, n_1, m_1, m_2, \theta),  
\edd
which gives us \eqref{a1pa4_chern_eqn}.  \qed \\

\begin{prp}
\label{a1pa5_chern}
The number $\N(\A_1\mc{P}\A_5, n_1, m_1, m_2, \theta)$ is given by 
\begin{align}
\N(\A_1\mc{P}\A_5, n_1, m_1, m_2, \theta) 
&= \N(\A_1\mc{P}\A_4, n_1, m_1, m_2, \theta+1) 
+4\N(\A_1\mc{P}\A_4, n_1, m_1+1, m_2, \theta) \nonumber \\  
& ~~ + 3 \N(\A_1\mc{P}\A_4, n_1, m_1, m_2, \theta)+3 \N(\A_1\mc{P}\A_4, n_1+1, m_1, m_2, \theta) \nonumber  \\  
& ~~ -2\N(\A_1\mc{P}\D_5, n_1, m_1, m_2, \theta) \nonumber \\ 
& ~~ -2 \N(\mc{P}\A_6, n_1, m_1, m_2, \theta)
-5 \N(\mc{P}\E_6, n_1, m_1, m_2, \theta).  \label{a1pa5_chern_eqn}
\end{align}
\end{prp}
\pf Let $\mu$ be a generic pseudocycle 
representing the homology class 
Poincar\'{e} dual to 
\[ c_1^{n_1} x_1^{m_1} x_2^{m_2} \lambda^{\theta} y^{\delta_L-(n_1+m_1+2m_2 + \theta+6)}.\]
We show in \cite{BM13_2pt_published}, that
\begin{align*}
\ov{\ov{\A}_1 \circ  \PP \A}_4 
&= \ov{\A}_1 \circ  \PP \A_4 \du \ov{\A}_1 \circ (\ov{\PP \A}_5 \cup \ov{\PP \D}_5) \du 
\Big( \Delta \ov{\PP \A}_6 \cup  
\Delta \ov{\mc{P}\mc{D}^{s}_7} \cup \Delta \ov{\PP \E}_6  \Big),  
\end{align*}
where 
\begin{align*}
\Delta \mc{P}\mc{D}^{s}_7 &:=  \{ ([f], q, l_q) \in \ov{\ov{\A}_1 \circ \PP \A}_4: 
\pi_2^*\us_{\PP \D_4}([f], q, l_q) =0, 
~\pi_2^*\us_{\PP \E_6}([f], q, l_q) \neq 0  \}.
\end{align*}
If $L$ is sufficiently $8$-ample, then 
the section 
$$ \pi_2^*\us_{\PP \A_5} : \ov{\ov{\A}_1 \circ \PP \A}_4 \cap \mu \lra \pi_2^*\UL_{\PP \A_5}$$
vanishes transversely on $\A_1 \circ \PP \A_5 \cap \mu$. 
We also 
show that this section  
vanishes on $\A_1 \circ \PP \D_5 \cap \mu$, $\Delta\PP \A_6 \cap \mu$ and $\Delta\PP \E_6 \cap \mu$ 
with a multiplicity of $2$, $2$ and $5$ respectively. 
Since the dimension of $\PP \D_7$ is one less than the 
dimension of $[\mu]$ and $\mu$ is generic, 
$\Delta \ov{\PP \D}_7 \cap \mu$ is empty. 
We also show in \cite{BM13_2pt_published} that 
$\Delta \ov{\mc{P}\mc{D}^{s}_7}$ is a subset of $\Delta \ov{\PP \D}_7$. 
Hence $\Delta \ov{\mc{P}\mc{D}^{s}_7} \cap \mu$ is also empty. 
Hence
\begin{align*}
\Big\langle e(\pi_2^*\UL_{\PP \A_5}), ~[\ov{\ov{\A}_1 \circ \PP \A}_4] \cap [\mu] \Big\rangle 
& = \Num(\A_1\PP\A_5,n_1,m_1, m_2, \theta) +2\Num(\A_1\PP\D_5,n_1,m_1, m_2, \theta)\\ 
& ~~ + 2 \Num(\PP \A_6, n_1, m_1, m_2, \theta)+ 5\Num(\PP \E_6, n_1 ,m_1, m_2, \theta),  
\end{align*}
which gives us \eqref{a1pa5_chern_eqn}.  \qed 

\begin{prp}
\label{a1pa6_chern}
The number $\N(\A_1\mc{P}\A_6, n_1, m_1, m_2, \theta)$ is given by 
\begin{align}
\N(\A_1\mc{P}\A_6, n_1, m_1, m_2, \theta)  
&=  
+6\N(\A_1\mc{P}\A_5, n_1, m_1+1, m_2, \theta) \nonumber \\  
& ~~ + 4 \N(\A_1\mc{P}\A_5, n_1, m_1, m_2, \theta)+4 \N(\A_1\mc{P}\A_5, n_1+1, m_1, m_2, \theta) \nonumber  \\  
& ~~ -4\N(\A_1\mc{P}\D_6, n_1, m_1, m_2, \theta) 
-3\N(\A_1\mc{P}\E_6, n_1, m_1, m_2, \theta)  \nonumber \\ 
& ~~ -2 \N(\mc{P}\A_7, n_1, m_1, m_2, \theta)
-6 \N(\mc{P}\E_7, n_1, m_1, m_2, \theta).  \label{a1pa6_chern_eqn}
\end{align}
\end{prp}
\pf Let $\mu$ be a generic pseudocycle 
representing the homology class Poincar\'{e} dual to 
\[ c_1^{n_1} x_1^{m_1} x_2^{m_2} \lambda^{\theta} y^{\delta_L-(n_1+m_1+2m_2 + \theta+5)}.\]
We show in \cite{BM13_2pt_published}, that
\begin{align*}
\ov{\ov{\A}_1 \circ  \PP \A}_5 & = 
\ov{\A}_1 \circ  \PP \A_5 \du \ov{\A}_1 \circ (\ov{\PP \A}_6 \cup \ov{\PP \D}_6 \cup \ov{\PP \E}_6 ) \du 
\Big( \Delta \ov{\PP \A}_7 \cup \Delta  \ov{\mq} \cup \Delta \ov{\PP \E}_7 \Big),
\end{align*}
where 
\begin{align*}
\Delta \mc{P}\mc{D}^{s}_8 &:=  \{ ([f], q, l_q) \in \ov{\ov{\A}_1 \circ \PP \A}_5: 
\pi_2^*\us_{\PP \D_4}([f], q, l_q) =0, 
~\pi_2^*\us_{\PP \E_6}([f], q, l_q) \neq 0  \}.
\end{align*}
If $L$ is sufficiently $9$-ample, then 
the section 
$$ \pi_2^*\us_{\PP \A_6}: 
\ov{\ov{\A}_1 \circ \PP \A}_5 \cap \mu \lra \pi_2^*\UL_{\PP \A_6}$$
vanishes transversely on $\A_1 \circ \PP \A_6$. 
We also show that 
this section  
vanishes on $\A_1 \circ \PP \D_6 \cap \mu$, $\A_1 \circ \PP \E_6 \cap \mu$ 
$\Delta\PP \A_7 \cap \mu$ and $\Delta\PP \E_7 \cap \mu$
with a multiplicity of $4$, $3$, $2$ and $6$ respectively. 
Since the dimension of $\PP \D_8$ is one less than the 
dimension of $[\mu]$ and $\mu$ is generic, 
$\Delta \ov{\PP \D}_8 \cap \mu$ is empty. 
We also show in \cite{BM13_2pt_published} that 
$\Delta \ov{\mc{P}\mc{D}^{s}_8}$ is a subset of $\Delta \ov{\PP \D}_8$. 
Hence $\Delta \ov{\mc{P}\mc{D}^{s}_8} \cap \mu$ is also empty. 
Hence 
\begin{align*}
\Big\langle e(\pi_2^*\UL_{\PP \A_6}), ~[\ov{\ov{\A}_1 \circ \PP \A}_5] \Big\rangle & = 
\Num(\A_1\PP\A_6,n_1,m_1, m_2, \theta) +4\Num(\A_1\PP\D_6,n_1,m_1, m_2, \theta)\\ 
& ~~ + 3\N(\A_1\mc{P}\E_6, n_1, m_1, m_2, \theta) + 2 \Num(\PP \A_7, n_1, m_1, m_2, \theta)\\
& ~~ + 6\Num(\PP \E_7, n_1 ,m_1, m_2, \theta)  
\end{align*}
which gives us \eqref{a1pa6_chern_eqn}.  \qed 

\begin{prp}
\label{a1pd4_chern}
The number $\N(\A_1\mc{P}\D_4, n_1, m_1, m_2, 0)$ is given by 
\begin{align}
\N(\A_1\mc{P}\D_4, n_1, m_1, m_2, 0) &= 2 \N(\A_1\mc{P}\A_3, n_1, m_1+1, m_2, 0) 
-2 \N(\A_1\mc{P}\A_3, n_1, m_1, m_2, 1) \nonumber \\ 
& +  \N(\A_1\mc{P}\A_3, n_1, m_1, m_2, 0) 
+\N(\A_1\mc{P}\A_3, n_1+1, m_1, m_2, 0) \nonumber \\ 
& -  2\N(\mc{P}\D_5, n_1, m_1, m_2, 0).\label{a1pd4_chern_eqn}
\end{align}
\end{prp}
\pf Let $\mu$ be a generic pseudocycle 
representing the homology class 
Poincar\'{e} dual to 
\[c_1^{n_1} x_1^{m_1} x_2^{m_2} y^{\delta_L-(n_1+m_1+2m_2 +5)}.\]
We show in \cite{BM13_2pt_published}, that
\begin{align*}
\ov{\ov{\A}_1 \circ  \PP \A}_3 &= 
\ov{\A}_1 \circ  \PP \A_3 \du \ov{\A}_1 \circ (\ov{\PP \A}_4 \cup \ov{\PP \D}_4) \du 
\Big( \Delta \ov{\PP \A}_5 \cup \Delta \ov{\PP \D^{\vee}_5} \Big).
\end{align*}
If $L$ is sufficiently $5$-ample, then 
the section 
$$ \pi_2^*\us_{\PP \D_4} : \ov{\ov{\A}_1 \circ \PP \A}_3 \cap \mu \lra \pi_2^*\UL_{\PP \D_4} $$
vanishes transversely on $\A_1 \circ \PP \D_4 \cap \mu$. 
It is easy to see that the section does not vanish on $\A_1 \circ \PP \A_4$ and $\Delta \PP \A_5$. We also show that 
the contribution of the section to the 
Euler class from the points of $\Delta \PP \D_5^{\vee} \cap \mu$ is $2$. 
Hence 
\begin{align}
\Big\langle e(\pi_2^*\UL_{\PP \D_4}), 
~[\ov{\ov{\A}_1 \circ \PP \A}_3] \cap [\mu] \Big\rangle & = \Num(\A_1\PP\D_4,n_1, m_1, m_2,0) + 
 2\Num(\D_5, n_1, m_1, m_2,0),\label{pd5_dual_prelim_a1_pa4}
\end{align}
giving us \eqref{a1pd4_chern_eqn}. \qed 

\begin{prp}
\label{a1pd4_theta_chern}
The number $\N(\A_1\mc{P}\D_4, n_1, m_1, m_2, 1)$ is given by 
\begin{align}
\N(\A_1\mc{P}\D_4, n_1, m_1, m_2, 1) &= 
\frac{1}{3} \N(\A_1 \mc{P}\D_4, n_1, m_1, m_2, 0) 
+\N(\A_1\mc{P}\D_4, n_1, m_1+1, m_2, 0) \nonumber \\ 
& + \frac{1}{3}\N(\A_1\mc{P}\D_4, n_1+1, m_1, m_2, 0).\label{a1pd4_theta_chern_eqn}
\end{align}
\end{prp}
\pf Let $\mu$ be a generic pseudocycle 
representing the homology class 
Poincar\'{e} dual to 
\[c_1^{n_1} x_1^{m_1} x_2^{m_2} \lambda y^{\delta_L-(n_1+m_1+2m_2 + 6)}. \]
We show in \cite{BM13_2pt_published}, that
\begin{align}
\ov{\ov{\A}_1 \circ  \hat{\D}_4} 
&= \ov{\A}_1 \circ  \hat{\D}_4^{\#} \du \ov{\A}_1 \circ \ov{\PP \D}_4 \du 
\Big( \Delta \ov{\hat{\D}^{\# \flat}_6} \Big). \nonumber
\end{align}
If $L$ is sufficiently $5$-ample, then 
the section 
$$ \pi_2^*\us_{\PP \A_3} : \ov{\ov{\A}_1 \circ  \hat{\D}_4} \cap \mu 
\lra \pi_2^*\UL_{\PP \A_3} $$
vanishes transversely on $\A_1 \circ \PP \D_4 \cap \mu$. 
By definition, the section does not vanish on $\Delta\hat{\D}^{\# \flat }_6 \cap \mu$. 
Hence 
\begin{align*}
\Big\langle e(\pi_2^*\UL_{\PP \A_3}), 
~[\ov{\ov{\A}_1 \circ  \hat{\D}^{\#}_4}] \cap [\mu] \Big\rangle & = \Num(\A_1\PP\D_4,n_1, m_1, m_2, 1)
\end{align*}
which gives us \eqref{a1pd4_theta_chern_eqn}. \qed

\begin{prp}
\label{a1pd4_theta_general_chern}
If $\theta >1$, the number $\N(\A_1\mc{P}\D_4, n_1, m_1, m_2, \theta)$ is given by 
\begin{align}
\N(\A_1\mc{P}\D_4, n_1, m_1, m_2, \theta) &= 
\N(\A_1 \mc{P}\D_4, n_1, m_1+1, m_2,\theta-1) \nonumber \\ 
& - \N(\A_1\mc{P}\D_4, n_1, m_1, m_2+1, \theta-2).  \label{a1pd4_theta_general_chern_eqn}
\end{align}
\end{prp}
\pf Follows immediately from the ring structure of $H^*(\mathbb{P} T\X)$. \qed 

\begin{prp}
\label{a1pd5_chern}
The number $\N(\A_1\mc{P}\D_5, n_1, m_1, m_2, \theta)$ is given by 
\begin{align}
\N(\A_1\mc{P}\D_5, n_1, m_1, m_2, \theta) &= \N(\A_1\mc{P}\D_4, n_1, m_1, m_2, \theta+1) 
+ \N(\A_1\mc{P}\D_4, n_1, m_1+1, m_2, \theta) \nonumber \\ 
& +\N(\A_1\mc{P}\D_4, n_1, m_1, m_2, \theta) + \N(\A_1\mc{P}\D_4, n_1+1, m_1, m_2, \theta) \nonumber \\
& -2 \N(\mc{P}\D_6, n_1, m_1, m_2, \theta).\label{a1pd5_chern_eqn}
\end{align}
\end{prp}
\pf Let $\mu$ be a generic pseudocycle 
representing the homology class 
Poincar\'{e} dual to 
\[c_1^{n_1} x_1^{m_1} x_2^{m_2} \lambda^{\theta} y^{\delta_L-(n_1+m_1+2m_2 + \theta+6)}.\]
We show in \cite{BM13_2pt_published}, that
\begin{align*}
\ov{\ov{\A}_1 \circ  \PP \D}_4 &= \ov{\A}_1 \circ  \PP \D_4 
\du \ov{\A}_1 \circ (\ov{\PP \D}_5 \cup \ov{\PP \D^{\vee}_5})  \du 
\Big( \Delta \ov{\PP \D}_6 \cup \Delta \ov{\mr} \Big).
\end{align*}
If $L$ is sufficiently $6$-ample, then 
the section 
$$ \pi_2^*\us_{\PP \D_5}^{\mathbb{L}} : \ov{\ov{\A}_1 \circ \PP \D}_4 \cap \mu 
\lra \pi_2^*\UL_{\PP \D_5} $$
vanishes transversely on $\A_1 \circ \PP \D_5$. 
Moreover, it does not vanish on $\ov{\A}_1 \circ \PP \D_5^{\vee} \cap \mu$ by definition. 
We also show that the  
the contribution of the section to the Euler class from the points of 
$\Delta\PP \D_6 \cap \mu$ is $2$.  
The section does not vanish on $\Delta \mr$ by definition (cf. \cite{BM13_2pt_published}). 
Since $\mu$ is generic, 
the section does not vanish on $\Delta \ov{\mr} \cap \mu$.  
Hence 
\begin{align*}
\Big\langle e(\pi_2^*\UL_{\PP \D_5}), 
~[\ov{\ov{\A}_1 \circ \PP \D}_4] \cap [\mu] \Big\rangle & = \Num(\A_1\PP\D_5,n_1,m_1, m_2, \theta) + 
2 \Num(\PP \D_6, n_1, m_1, m_2, \theta)
\end{align*}
which gives us \eqref{a1pd5_chern_eqn}. \qed 

\begin{prp}
\label{a1pd6_chern}
The numbers $\N(\A_1\mc{P}\D_6, n_1, m_1, m_2, \theta)$ and 
$\N(\A_1\mc{P}\E_6, n_1, m_1, m_2, \theta)$
are 
\begin{align}
\N(\A_1\mc{P}\D_6, n_1, m_1, m_2, \theta) &= 4 \N(\A_1\mc{P}\D_5, n_1, m_1, m_2, \theta+1) 
+ \N(\A_1\mc{P}\D_5, n_1, m_1, m_2, \theta) \nonumber \\ 
& +\N(\A_1\mc{P}\D_5, n_1+1, m_1, m_2, \theta)  \nonumber \\
& -2 \N(\mc{P}\D_7, n_1, m_1, m_2, \theta)-\N(\mc{P}\E_7, n_1, m_1, m_2, \theta),  \label{a1pd6_chern_eqn} 
\end{align}
\begin{align}
\N(\A_1\mc{P}\E_6, n_1, m_1, m_2, \theta) &= 2 \N(\A_1\mc{P}\D_5, n_1, m_1+1, m_2, \theta) 
-\N(\A_1\mc{P}\D_5, n_1, m_1, m_2, \theta+1) \nonumber \\ 
& +\N(\A_1\mc{P}\D_5, n_1, m_1, m_2, \theta) +\N(\A_1\mc{P}\D_5, n_1+1, m_1, m_2, \theta)  \nonumber \\
& -\N(\mc{P}\E_7, n_1, m_1, m_2, \theta). \label{a1pe6_chern_eqn} 
\end{align}
\end{prp}
\pf Let $\mu$ be a generic pseudocycle 
representing the homology class 
Poincar\'{e} dual to 
\[ c_1^{n_1} x_1^{m_1} x_2^{m_2} \lambda^{\theta} y^{\delta_L-(n_1+m_1+2m_2 + \theta+7)}. \]
We show in \cite{BM13_2pt_published}, that
\begin{align*}
\ov{\ov{\A}_1 \circ  \PP \D}_5 &= \ov{\A}_1 \circ  \PP \D_5 
\du \ov{\A}_1 \circ (\ov{\PP \D}_6 \cup \ov{\PP \E}_6) \du 
\Big( \Delta \ov{\PP \D}_7 \cup \Delta \ov{\PP \E}_7 \Big). 
\end{align*}
If $L$ is sufficiently $7$-ample or $6$-ample, then 
the sections 
\[ \pi_2^*\us_{\PP \D_6} : \ov{\ov{\A}_1 \circ \PP \D}_5 \cap \mu \lra \pi_2^*\UL_{\PP \D_6}, 
~~\pi_2^*\us_{\PP \E_6}  : \ov{\ov{\A}_1 \circ \PP \D}_5 \cap \mu \lra \pi_2^*\UL_{\PP \E_6} \]
vanish transversely on $\A_1 \circ \PP \D_6 \cap \mu$ and $\A_1 \circ \PP \E_6 \cap \mu$ 
respectively. 
Moreover, they do not vanish on $\A_1 \circ \PP \E_6 \cap \mu$ and $\A_1 \circ \PP \D_6 \cap \mu$ respectively. 
We also show that 
the contribution of the sections  $\pi_2^*\us_{\PP \D_6}$ 
to the Euler class from the points of $\Delta\PP \D_7 \cap \mu$ and $\Delta \PP \E_7 \cap \mu$ are 
$2$ and $1$ respectively. 
We also show that the contribution of the section    
$\pi_2^*\us_{\PP \E_6}$ from the points of $\Delta \PP \E_7 \cap \mu$ is $1$; 
moreover it does not vanish on $\Delta \PP \D_7 \cap \mu$. 
Hence 
\begin{align*}
\Big\langle e(\pi_2^*\UL_{\PP \D_6} ), ~[\ov{\ov{\A}_1 \circ \PP \D}_5] \cap \mu \Big\rangle & = 
\Num(\A_1\PP\D_6,n_1,m_1, m_2, \theta) \\ 
& + 2 \Num(\PP \D_7, n_1, m_1, m_2, \theta)+\Num(\PP \E_7, n_1, m_1, m_2, \theta), \\
\Big\langle e(\pi_2^*\UL_{\PP \E_6} ), ~[\ov{\ov{\A}_1 \circ \PP \D}_5] \cap \mu \Big\rangle & = 
\Num(\A_1\PP\E_6,n_1,m_1, m_2, \theta) + \Num(\PP \E_7, n_1, m_1, m_2, \theta),
\end{align*}
which give us \eqref{a1pd6_chern_eqn} and \eqref{a1pe6_chern_eqn} respectively. \qed


\bibliographystyle{siam}
\bibliography{Myref_bib.bib}

\vspace*{0.4cm}

\hf {\small D}{\scriptsize EPARTMENT OF }{\small M}{\scriptsize ATHEMATICS, }{\small RKM V}{\scriptsize IVEKANANDA }{\small U}{\scriptsize NIVERSITY, }{\small H}{\scriptsize OWRAH, }{\small WB} {\footnotesize 711202, }{\small INDIA}\\
\hf{\it E-mail address} : \texttt{somnath@maths.rkmvu.ac.in}\\

\hf {\small D}{\scriptsize EPARTMENT OF }{\small M}{\scriptsize ATHEMATICS,}
{\small T}{\scriptsize ATA}
{\small I}{\scriptsize NSTITUTE OF}
{\small F}{\scriptsize UNDAMENTAL}
{\small R}{\scriptsize ESEARCH,} 
{\small M}{\scriptsize UMBAI}
{\footnotesize 400005, }{\small INDIA}\\
\hf{\it E-mail address} : \texttt{ritwikm@math.tifr.res.in}\\[0.2cm]

\end{document}

%% file: New_Macros_by_Ritwik_and_Somnath.tex
\def \hf{\hspace*{0.5cm}}                      
\def\bge{\begin{equation}}                
\def\ede{\end{equation}}                
\def\bgd{\begin{displaymath}}         
\def\edd{\end{displaymath}}            
\def\bgee{\begin{equation*}}           
\def\edee{\end{equation*}}           

\def \ni{\noindent}
\def \hf{\hspace*{0.5cm}}                      
\def\lra{\longrightarrow}

\def\lan{\langle}
\def\ran{\rangle}

\def \C{\mathbb{C}}

\def \A{\mathcal{A}}
\def \D{\mathcal{D}}
\def \ov{\overline}

\def \P{\mathbb{P}}
\def \N{\mathcal{N}}
\def \A{\mathcal{A}}
\def \E{\mathcal{E}}

\def \X{X}

\def \pf{\noindent \textbf{Proof:  }}

\def \mf{\mathfrak}
\def \mc{\mathcal}

\providecommand \hf{\hspace*{0.5cm}}                      
\providecommand\bge{\begin{equation}}                
\providecommand\ede{\end{equation}}                
\providecommand\bgd{\begin{displaymath}}         
\providecommand\edd{\end{displaymath}}            
\providecommand\bgee{\begin{equation*}}           
\providecommand\edee{\end{equation*}}           

\providecommand \ni{\noindent}

\providecommand\lra{\longrightarrow}

\providecommand\lan{\langle}
\providecommand\ran{\rangle}

\providecommand\BA{\begin{eqnarray}}
\providecommand\EA{\end{eqnarray}}
\providecommand\BAA{\begin{eqnarray*}}
\providecommand\EAA{\end{eqnarray*}}
\providecommand\Bal{\begin{align*}}
\providecommand\Eal{\end{align*}}

\providecommand \C{\mathbb{C}}

\providecommand \P{\mathbb{P}}
\providecommand \A{\mathcal{A}}
\providecommand \E{\mathcal{E}}
\providecommand \Num{\mathcal{N}}
\providecommand \PP{\mathcal{P}}
\providecommand \f{\frac}

\providecommand \UL{\mathbb{L}} 
\providecommand \DL{\mathcal{L}}
\providecommand \UV{\mathbb{V}}
\providecommand \DV{\mathcal{V}}

\providecommand \ds{\psi}
\providecommand \us{\Psi}

\providecommand \ct{\mathcal{C}}

\providecommand \D{\mathcal{D}}

\providecommand \X{\mathfrak{X}}

\providecommand \pf{\noindent \textbf{Proof:  }}
\providecommand \ov{\overline}

\providecommand \XC{\mathcal{X}}

\providecommand \U{\mathcal{U}}

\providecommand \N{\nabla}

\providecommand \du{\sqcup}

\providecommand \mq{\PP \D_8^{s}}
\providecommand \mr{\PP \D_6^{\vee s}}